\documentclass[a4paper]{amsart}

\usepackage[pdftex]{graphicx}
\usepackage{amsmath,amssymb,amsthm}
\usepackage{mathrsfs}
\usepackage{mathtools}
\usepackage{xcolor}
\usepackage[T1]{fontenc}
\usepackage{tikz}
\usetikzlibrary{arrows.meta}
\usepackage{hyperref}

\theoremstyle{plain}
\newtheorem{thm}{Theorem}[section]
\newtheorem{prop}[thm]{Proposition}
\newtheorem{lem}[thm]{Lemma}
\newtheorem{cor}[thm]{Corollary}

\theoremstyle{definition}
\newtheorem{dfn}[thm]{Definition}

\theoremstyle{remark}
\newtheorem{rem}[thm]{Remark}

\renewcommand{\theenumi}{\roman{enumi}}
\newcommand{\C}{\mathbb{C}}
\newcommand{\R}{\mathbb{R}}
\newcommand{\Z}{\mathbb{Z}}
\newcommand{\Q}{\mathbb{Q}}
\newcommand{\A}{\mathbb{A}}
\newcommand{\Span}{\operatorname{Span}}
\newcommand{\Spec}{\operatorname{Spec}}
\newcommand{\FB}{\operatorname{FB}}
\newcommand{\wt}{\operatorname{wt}}
\newcommand{\Hom}{\operatorname{Hom}}
\newcommand{\Hilb}{\operatorname{Hilb}}
\newcommand{\tor}{\operatorname{tor}}
\newcommand{\GL}{\operatorname{GL}}
\newcommand{\SL}{\operatorname{SL}}
\newcommand{\diag}{\operatorname{diag}}

\begin{document}
\title{$G$-Hilbert schemes of canonical abelian quotient threefolds}
\author{Yutaro Kaijima}
\address{Department of Mathematics, Graduate School of Science, the University of Osaka, Toyonaka, Osaka 560-0043, Japan}
\email{u797381f@ecs.osaka-u.ac.jp}
\date{}

\begin{abstract}
Ishida and Iwashita showed that the canonical abelian quotient singularities of dimension three that are not Gorenstein fall into four families.
For the two families that have not been treated before, we extend a method of K\k{e}dzierski and compute the toric charts of the $G$-Hilbert scheme explicitly.
As a result, we show that the $G$-Hilbert scheme has at worst conifold singularities.
Combined with known results, this shows that for every three-dimensional canonical abelian quotient singularity, the $G$-Hilbert scheme has at worst conifold singularities, every essential divisor over the singularity appears on it as a prime divisor, and its limit F-blowup is a resolution of the singularity.
These results also hold for every three-dimensional $\Q$-factorial canonical toric singularity in any characteristic, if we replace the $G$-Hilbert scheme by the limit F-blowup of the singularity.
\end{abstract}

\maketitle

\section{Introduction}

Let $k$ be an algebraically closed field, let $G$ be a finite abelian subgroup of $\GL(3,k)$ of order prime to the characteristic of $k$, and let $X=\A_k^3/G$ be the associated quotient variety.
Ishida and Iwashita \cite{II86} classified the groups $G$ for which $X$ has at worst a canonical singularity.
Apart from the case $G\subset\SL(3,k)$, there are exactly four families of groups.
We call them types (2)--(5), as in Proposition \ref{canonical classification}.
Here, for a primitive $r$-th root of unity $\varepsilon\in k$, we say that $G$ is of type $\frac1r(a_1,a_2,a_3)$ if $G$ is the cyclic group generated by $\diag(\varepsilon^{a_1},\varepsilon^{a_2},\varepsilon^{a_3})$.
Two of these families have already been studied in detail.
Type (2) is the terminal case and consists of the groups of type $\frac1r(1,a,r-a)$ with $\gcd(a,r)=1$ (see \cite{MS84}).
Type (5) consists of the two sporadic groups of types $\frac1{9}(1,4,7)$ and $\frac1{14}(1,9,11)$.
This paper treats the remaining two families.
Type (3) consists of the groups of type $\frac1r(1,r-1,a)$ with $\gcd(a,r)>1$, and type (4) consists of the groups of type $\frac1{4t}(1,2t+1,4t-2)$ with $t\geq2$.

The $G$-Hilbert scheme $\Hilb^G(\A_k^3)$ is the distinguished component of the moduli space of $G$-clusters \cite{IN96,Nak01}.
It is a natural candidate for a resolution of $X$.
When $G\subset\SL(3,k)$, it is a crepant resolution of $X$. This was shown by Nakamura \cite{Nak01} for abelian $G$, and by Bridgeland, King and Reid \cite{BKR01} for arbitrary finite $G$.
For a finite abelian subgroup $G\subset\GL(d,k)$, $\Hilb^G(\A_k^d)$ is a toric variety.
Its toric structure has been described explicitly in several other cases \cite{Kid01,CR02,Seb07}, and in general by Craw, Maclagan and Thomas using initial ideals \cite{CMT07a,CMT07b}.

For type (2), K\k{e}dzierski \cite{Ked11} showed that $\Hilb^G(\A_k^3)$ has at worst conifold singularities.
The proof describes the toric charts explicitly in terms of \emph{$G$-graphs}, which are classified by their number of \emph{valleys}.
For type (5), Sato and Sato \cite{SS23} showed that the normalization of $\Hilb^G(\A_k^3)$ is a resolution of $X$, and we check in Remark \ref{type5 smooth} that $\Hilb^G(\A_k^3)$ itself is smooth.

We follow K\k{e}dzierski's method.
For the groups of types (3) and (4), we classify the $G$-graphs by their number of valleys.
For each $G$-graph $\Gamma$ with $\dim\sigma(\Gamma)=3$, we then compute the coordinate ring $k[S(\Gamma)]$ of the corresponding toric chart of $\Hilb^G(\A_k^3)$. These charts cover $\Hilb^G(\A_k^3)$.
In every case, $k[S(\Gamma)]$ is isomorphic to $k[x,y,z]$ or to $k[x,y,z,w]/(xy-zw)$.
From this, we obtain our first main result.

\begin{thm}[Corollaries \ref{Hilb conifold} and \ref{Hilb conifold section4}]\label{thm:intro conifold}
Let $G\subset\GL(3,k)$ be of type (3) or (4).
Then $\Hilb^G(\A_k^3)$ has at worst conifold singularities.
In particular, $\Hilb^G(\A_k^3)$ is normal.
\end{thm}

Together with the known cases, this gives a description of $\Hilb^G(\A_k^3)$ for all three-dimensional canonical abelian quotient singularities.
For a toric variety $Y$ and a positive integer $l$, let $\FB_{(l)}(Y)$ denote the $(l)$-th F-blowup of $Y$, that is, the universal flattening of a Frobenius-like morphism \cite{Yas09,Yas12}.
Let $\FB_{(\infty)}(Y)$ denote its stabilization for $l\gg0$, which we call the limit F-blowup of $Y$. See \cite{Har12,HSY13,Har15} for surfaces and \cite{TY09,CMDY25,KY26,CMKY26} for higher-dimensional toric singularities.
We combine Theorem \ref{thm:intro conifold} with three results: the theory of essential divisors of toric singularities by Bouvier and Gonzalez-Sprinberg \cite{BGS95}, its refinement by Ch\'avez-Mart\'inez, Duarte and Yasuda \cite{CMDY25}, and the theory of F-blowups by Yasuda \cite{Yas12}.
As a result, we obtain our second main result.

\begin{thm}[Theorem \ref{Hilb canonical} and Corollary \ref{FB FB resolution}]\label{thm:intro canonical}
Let $G\subset\GL(3,k)$ be a finite abelian subgroup containing no pseudo-reflections, of order prime to the characteristic of $k$, such that $X=\A_k^3/G$ has at worst a canonical singularity.
Then:
\begin{enumerate}
\item $\Hilb^G(\A_k^3)$ has at worst conifold singularities; in particular, it is normal;
\item every essential divisor over $X$ appears on $\Hilb^G(\A_k^3)$ as a prime divisor;
\item $\FB_{(l)}(\Hilb^G(\A_k^3))$ is a resolution of $X$ for every $l\geq2$.
\end{enumerate}
\end{thm}

Essential divisors come from the work of Nash \cite{Nash95} on arc spaces.
Every irreducible component of the space of arcs through the singular locus determines an essential divisor.
The Nash problem asks whether every essential divisor arises in this way.
The answer is affirmative for toric varieties \cite{IK03} and for surfaces \cite{FdBPP12}, but negative in general in dimension at least three \cite{IK03,dF13,JK13} (see \cite{PS15} for a survey).
Thus essential divisors form an important class of divisors over a singularity.
It is natural to ask which partial resolutions contain all of them.
Assertion (ii) gives an affirmative answer to \cite[Problem 1.1]{CMDY25} for these singularities.
This was known for $G\subset\SL(3,k)$ \cite[Corollary 6.13]{CMDY25} and for type (2) \cite[Theorem 5.12]{CMKY26}. In general, it was only known that the centers of essential divisors on $\Hilb^G(\A_k^3)$ have positive dimension \cite[Corollary 5.18]{CMKY26}. Moreover, the method of \cite{CMDY25} does not apply, since $X$ has no moderate toric resolution unless $X$ is Gorenstein \cite[Corollary 3.15]{KY26}.

By \cite{Yas12,TY09}, we have $\Hilb^G(\A_k^3)\cong\FB_{(\infty)}(X)$.
Hence assertion (iii) says that $\FB_{(\infty)}(\FB_{(\infty)}(X))$ is a resolution of $X$.
This gives a positive answer to \cite[Question 1.4]{Yas12} for $X$ (compare the two-step construction of \cite{IIN13}).
The key point is that for $l\geq2$, the F-blowup $\FB_{(l)}$ of a conifold singularity is its blowup at the singular point.
The analogue of Theorem \ref{thm:intro canonical} holds for $\FB_{(\infty)}(X)$ for every three-dimensional $\Q$-factorial canonical toric singularity $X$, in any characteristic (Remark \ref{FB conifold}).

The paper is organized as follows.
In Section \ref{sec:prelim}, we recall Nakamura's description of $\Hilb^G(\A_k^d)$ in terms of $G$-graphs, the classification of Ishida and Iwashita, and the notion of essential divisors.
In Sections \ref{sec:type3} and \ref{sec:type4}, we classify the $G$-graphs for the groups of types (3) and (4), respectively, compute the corresponding charts, and prove Theorem \ref{thm:intro conifold}.
In Section \ref{sec:FB}, we study the F-blowups of the conifold singularity, and prove Theorem \ref{thm:intro canonical} by combining these results with the known cases.

\section*{Acknowledgments}

The author would like to thank Yusuke Sato, who suggested the original idea underlying this paper, for invaluable guidance, encouragement, and many helpful discussions.
The author would also like to thank Takehiko Yasuda, the author's advisor, for insightful comments and valuable advice, and for reading through this manuscript.
The author is grateful to Linghu Fan and Yukari Ito for many helpful discussions as well.
This work was supported by JST SPRING, Grant Number JPMJSP2138.

\section*{Declaration of AI Use}

The author used Claude (Anthropic) to assist with mathematical discussions and the preparation of this manuscript.
All mathematical arguments were independently verified by the author, who takes full responsibility for the content of this paper.

\section{Preliminaries}\label{sec:prelim}

In this section, we fix notation for abelian quotient singularities and recall Nakamura's description of $\Hilb^G(\A_k^d)$ by $G$-graphs.
We also give some lemmas used in Sections \ref{sec:type3} and \ref{sec:type4}.
Finally, we recall the classification of canonical abelian quotient threefolds and the notion of essential divisors.

\subsection{\texorpdfstring{Abelian quotient singularities (See \cite[Section 1.1]{Nak01})}{Abelian quotient singularities}}

Let $k$ be an algebraically closed field.
Let $G$ be a finite abelian subgroup of $\GL(d,k)$ of order $r$ prime to the characteristic of $k$, and fix a primitive $r$-th root of unity $\varepsilon\in k$.
Without loss of generality, we may assume that every element of $G$ is a diagonal matrix.

Let $N_0$ be a free $\Z$-module with basis $e_1,\ldots,e_d$, let $M_0\coloneqq\Hom(N_0,\Z)$ be its dual with dual basis $f_1,\ldots,f_d$, and let $\Delta\coloneqq\sum_{i=1}^d\R_{\geq0}e_i$.
Given an element $g$ of $G$, we write $g=\diag(\varepsilon^{a_1},\ldots,\varepsilon^{a_d})$ with $0\leq a_i<r$, and set $\bar g\coloneqq(a_1/r)e_1+\cdots+(a_d/r)e_d\in N_0\otimes_\Z\Q$.
We define
\[
N\coloneqq N_0+\sum_{g\in G}\Z\bar g,\quad M\coloneqq\Hom_\Z(N,\Z).
\]
Note that $N_0\otimes_\Z\R=N\otimes_\Z\R$ and $M_0\otimes_\Z\R=M\otimes_\Z\R$, and we denote these spaces by $N_\R$ and $M_\R$, respectively.

Let $G^\vee$ be the character group of $G$.
The group $G$ acts on regular functions on the affine space $\A_{k}^d$ on the left, by setting $(g\cdot f)(P)=f(g^{-1}\cdot P)$, where $g\in G$, $P\in\A_{k}^d$ and $f$ is a regular function on $\A_{k}^d$.
For $m=\sum_{i=1}^d m_if_i\in M_0$, write
\[
x^m=x_1^{m_1}\cdots x_d^{m_d}.
\]
The action of $G$ on Laurent monomials induces a homomorphism
\[
\wt\colon M_0\to G^\vee,
\]
defined by
\[
g\cdot x^m=\wt(m)(g)x^m.
\]
We identify $M_0$ with the group of Laurent monomials in $k[x_1^{\pm1},\ldots,x_d^{\pm1}]$ via $m\mapsto x^m$, write $\wt(x^m)\coloneqq\wt(m)$, and use the additive and multiplicative notations interchangeably.
Since $\langle\bar g,m\rangle\in\Z$ if and only if $\wt(m)(g)=1$, the lattice $M\subset M_0$ is the kernel of $\wt$, that is, the lattice of exponents of $G$-invariant Laurent monomials.
For a positive integer $r$ and integers $a_1,\ldots,a_d$, we say that $G$ is of \emph{type} $\frac{1}{r}(a_1,\ldots,a_d)$ if $G$ is the cyclic group generated by $\diag(\varepsilon^{a_1},\ldots,\varepsilon^{a_d})$.
When $d=3$, we write $x,y,z$ for $x_1,x_2,x_3$.

\subsection{\texorpdfstring{$G$-Hilbert schemes}{G-Hilbert schemes}}

\begin{dfn}[{\cite{IN96,Nak01,BKR01}}]
A \emph{$G$-cluster} is a $G$-invariant closed subscheme $Z\subset\A_k^d$ such that $H^0(\mathcal{O}_Z)$ is isomorphic to the regular representation of $G$ as a representation of $G$.
The \emph{$G$-Hilbert scheme} $G$-$\Hilb(\A_k^d)$ is the fine moduli space parametrizing $G$-clusters. It has a Hilbert--Chow morphism $\tau\colon G$-$\Hilb(\A_k^d)\to\A_k^d/G$, which sends a $G$-cluster to the orbit that supports it.
We write $\Hilb^G(\A_k^d)$ for the irreducible component of $G$-$\Hilb(\A_k^d)$ on which $\tau$ is birational, that is, the component containing the free $G$-orbits.
\end{dfn}

For abelian $G$, Nakamura \cite{Nak01} described an open covering of $\Hilb^G(\A_k^d)$ by affine toric varieties in terms of the following combinatorial objects.

Let $M_0^+\coloneqq\bigoplus_{i=1}^d\Z_{\ge0}f_i$.
We identify $M_0^+$ with the semigroup of exponents of monomials in $k[x_1,\ldots,x_d]$.

\begin{dfn}[{\cite[Definition 1.4]{Nak01}}]
A subset $\Gamma$ of monomials in $k[x_1,\ldots,x_d]$ is called a \emph{$G$-graph} if
\begin{enumerate}
\item it contains the constant monomial 1,
\item if $vw\in\Gamma$ then $v\in\Gamma$ and $w\in\Gamma$,
\item $\wt|_\Gamma\colon\Gamma\to G^\vee$ is a bijection.
\end{enumerate}
\end{dfn}

\begin{dfn}[{\cite[Definition 1.5]{Nak01}}]
Let $\Gamma$ be a $G$-graph.
We define the weight map $\wt_\Gamma\colon M_0^+\to\Gamma$ as follows.
For any $u\in M_0^+$ there exists a unique element $v\in\Gamma$ such that $\wt(v)=\wt(u)$.
We denote $v$ by $\wt_\Gamma(u)$ and also define $s_\Gamma(u)\coloneqq u\cdot\wt_\Gamma(u)^{-1}$.
Then we define cones $\sigma(\Gamma)$ and $\sigma(\Gamma)^\vee$ by
\begin{align*}
\sigma(\Gamma)&\coloneqq\{\alpha\in N_\R\mid\langle\alpha,s_\Gamma(v)\rangle\geq0,\forall v\in M_0^+\}, \\
\sigma(\Gamma)^\vee&\coloneqq\{v\in M_\R\mid\langle\alpha,v\rangle\geq0,\forall\alpha\in\sigma(\Gamma)\},
\end{align*}
where $\langle\cdot,\cdot\rangle$ denotes the pairing between $N_\R$ and $M_\R$.
Since $\wt(s_\Gamma(v))$ is trivial, $s_\Gamma(v)$ lies in $M$. We define $S(\Gamma)$ to be the subsemigroup of $M$ generated by the set $\{s_\Gamma(v)\mid v\in M_0^+\}$.
\end{dfn}

The following elementary observation will be used repeatedly.

\begin{lem}[{cf.\ \cite[Remark 3.1]{Ked11}}]\label{not divisible}
Let $\Gamma$ be a $G$-graph, let $u\in M_0^+\setminus\Gamma$, and let $1\leq i\leq d$ be such that $u/x_i\in\Gamma$.
Then $\wt_\Gamma(u)$ is not divisible by $x_i$.
\end{lem}

\begin{proof}
Suppose $\wt_\Gamma(u)=x_iw$ for some monomial $w$.
Since $\Gamma$ is closed under taking divisors, $w\in\Gamma$.
Since $\wt(u/x_i)=\wt(\wt_\Gamma(u)/x_i)=\wt(w)$ and $u/x_i,w\in\Gamma$, we have $u/x_i=w$.
Hence $\wt_\Gamma(u)=x_iw=u$, contradicting $u\notin\Gamma$.
\end{proof}

To compute $k[S(\Gamma)]$, we use the following finite sets of generators of $S(\Gamma)$.

\begin{prop}[{\cite[Lemma 1.8]{Nak01}}]\label{monoid generator}
Let $\Gamma$ be a $G$-graph with $\dim\sigma(\Gamma)=d$ and let $A\subset M_0^+\setminus\Gamma$ be a finite set such that $M_0^+\setminus\Gamma=A\cdot M_0^+$.
Then $S(\Gamma)$ is generated by the finite set $\{s_\Gamma(v)\mid v\in A\}$ as a semigroup.
\end{prop}

Without the assumption on $\dim\sigma(\Gamma)$, we still have the following finite set of generators.

\begin{lem}\label{generators B}
Let $\Gamma$ be a $G$-graph. Then $S(\Gamma)$ is generated as a semigroup by the finite set
\[
B_\Gamma\coloneqq\{s_\Gamma(x_i\cdot w)\mid w\in\Gamma,\ 1\leq i\leq d,\ x_i\cdot w\notin\Gamma\}.
\]
\end{lem}

\begin{proof}
For $u\in M_0^+$ and $1\leq i\leq d$, we have $\wt_\Gamma(x_i\cdot u)=\wt_\Gamma(x_i\cdot\wt_\Gamma(u))$, and hence
\[
s_\Gamma(x_i\cdot u)=s_\Gamma(u)\,s_\Gamma(x_i\cdot\wt_\Gamma(u)).
\]
Here $s_\Gamma(x_i\cdot\wt_\Gamma(u))$ is either $1$ (if $x_i\cdot\wt_\Gamma(u)\in\Gamma$) or an element of $B_\Gamma$.
By induction on the degree of $v\in M_0^+$, it follows that every $s_\Gamma(v)$ lies in the semigroup generated by $B_\Gamma$.
\end{proof}

Once a small set of generators of $S(\Gamma)$ is found, $k[S(\Gamma)]$ is determined by the following lemma.

\begin{lem}[{cf.\ \cite[Lemmas 1.7 and 3.7]{Nak01}}]\label{basis}
Let $\Gamma$ be a $G$-graph. Then $S(\Gamma)$ generates $M$ as a group.
Consequently, if $d=3$ and $S(\Gamma)$ is generated by three elements, then $k[S(\Gamma)]\cong k[x,y,z]$. If $d=3$ and $S(\Gamma)$ is generated by four elements $u_1,\ldots,u_4$ with $u_1u_2=u_3u_4$, then $k[S(\Gamma)]\cong k[x,y,z,w]/(xy-zw)$.
\end{lem}

\begin{proof}
The first assertion is \cite[Lemma 1.7]{Nak01}.
Hence three generators of $S(\Gamma)$ generate $M\cong\Z^3$, so they form a $\Z$-basis of $M$.
In the second case, the kernel of $k[t_1,\ldots,t_4]\to k[S(\Gamma)]$, $t_i\mapsto u_i$, is a prime ideal of height $4-3=1$.
It contains the prime element $t_1t_2-t_3t_4$, so it is equal to $(t_1t_2-t_3t_4)$.
\end{proof}

The semigroup rings $k[S(\Gamma)]$ give the charts of $\Hilb^G(\A_k^d)$, as follows.

\begin{dfn}[{\cite[Definition 1.9]{Nak01}}]
For any $G$-graph $\Gamma$, we define
\[
U_{\sigma(\Gamma)}\coloneqq\Spec k[\sigma(\Gamma)^\vee\cap M],\quad V(\Gamma)\coloneqq\Spec k[S(\Gamma)].
\]
We note that $U_{\sigma(\Gamma)}$ is the normalization of $V(\Gamma)$.
\end{dfn}

\begin{thm}[{\cite[Theorem 2.11]{Nak01}, \cite[Theorem 5.2]{CMT07b}}]\label{Nakamura fan}
Let $G$ be a finite abelian subgroup of $\GL(d,k)$.
When $\Gamma$ ranges over all $G$-graphs, the set of all faces of all $d$-dimensional cones $\sigma(\Gamma)$ forms a fan in $N_\R$ supported on $\Delta$.
Moreover, the toric variety defined by this fan is isomorphic to the normalization of $\Hilb^G(\A_k^d)$.
Furthermore, the toric varieties $V(\Gamma)$ form an open covering of $\Hilb^G(\A_k^d)$ when $\Gamma$ ranges over all $G$-graphs with $\dim\sigma(\Gamma)=d$.
\end{thm}

\begin{rem}\label{full-dim charts}
In \cite{Nak01}, Theorem 2.11 and Lemma 1.8 (Proposition \ref{monoid generator} above) are stated without the assumption $\dim\sigma(\Gamma)=d$.
Without this assumption, they are false in general (see \cite[Example 4.12]{CMT07b}).
Accordingly, we compute $k[S(\Gamma)]$ only for the $G$-graphs with $\dim\sigma(\Gamma)=3$, which suffice by the above theorem.
\end{rem}

Finally, we introduce a notation for describing $G$-graphs.

\begin{dfn}[See \cite{ISS26}]
We say that a $G$-graph $\Gamma$ is spanned by monomials $u_1,\ldots,u_n$ if $\Gamma$ consists of all monomials dividing at least one of $u_1,\ldots,u_n$, and we write $\Gamma=\Span(u_1,\ldots,u_n)$ (not to be confused with a linear span).
\end{dfn}

\subsection{Classification of canonical abelian quotient singularities}

We recall the classification of three-dimensional canonical abelian quotient singularities, which reduces our problem to finitely many families of groups.

\begin{prop}[{\cite[Theorem 4.1]{II86}}]\label{canonical classification}
Let $G$ be a finite cyclic subgroup of $\GL(3,\C)$ containing no pseudo-reflections. Then $\A_\C^3/G$ has at worst a canonical singularity if and only if, after changing the generator of $G$ and permuting the coordinates if necessary,
\begin{enumerate}
\renewcommand{\theenumi}{\arabic{enumi}}
\item $G$ is of type $\frac{1}{r}(a,b,c)$ with $a+b+c\equiv0\pmod{r}$,
\item $G$ is of type $\frac{1}{r}(1,a,r-a)$ with $\gcd(a,r)=1$ and $0<a<r$,
\item $G$ is of type $\frac{1}{r}(1,r-1,a)$ with $\gcd(a,r)>1$ and $0<a<r$,
\item $G$ is of type $\frac{1}{4t}(1,2t+1,4t-2)$ with $t\geq2$, or
\item $G$ is of type $\frac{1}{9}(1,4,7)$ or $\frac{1}{14}(1,9,11)$.
\end{enumerate}
\end{prop}

\begin{rem}
By \cite{II86}, the three-dimensional non-Gorenstein canonical abelian quotient singularities are exactly those of types (2)--(5). That is, for a finite abelian subgroup $G\subset\GL(3,\C)$ containing no pseudo-reflections, the quotient $\A_\C^3/G$ has at worst a canonical singularity if and only if $G\subset\SL(3,\C)$ or $G$ is of one of the types (2)--(5).
Since this classification is carried out in terms of the defining cone of $\A_k^3/G$ (see \cite[Theorem 3.6]{II86}), the same holds over any algebraically closed field $k$ of characteristic prime to the order of $G$.
\end{rem}

The case $G\subset\SL(3,k)$ and types (2) and (5) are known (see the introduction and Remark \ref{type5 smooth}). In this paper, we focus on types (3) and (4).

\begin{rem}\label{type5 smooth}
Let $G$ be of type (5), that is, of type $\frac19(1,4,7)$ or $\frac1{14}(1,9,11)$.
Then $\Hilb^G(\A_k^3)$ is smooth.
Indeed, we enumerated the $G$-graphs by computer and computed the set $B_\Gamma$ of Lemma \ref{generators B} for each of them.
For each $G$-graph $\Gamma$ with $\dim\sigma(\Gamma)=3$, we checked that $B_\Gamma$ contains a $\Z$-basis of $M$ and that every element of $B_\Gamma$ is a non-negative integral combination of this basis.
Hence $V(\Gamma)\cong\A_k^3$.
The program and its output are available as ancillary files accompanying the arXiv version of this paper.
\end{rem}

\subsection{Essential divisors}

We recall the notion of essential divisors, and show that every essential divisor over a $\Q$-factorial affine toric variety appears on any toric partial resolution with at worst conifold singularities (Lemma \ref{conifold ess}).

\begin{dfn}[{\cite[Definition 2.3]{IK03}}]
Let $X$ be a normal variety.
We say that an exceptional divisor $E$ over $X$ is \emph{essential} if for every resolution $f\colon Z\to X$, the center of $E$ on $Z$ is an irreducible component of $f^{-1}(X_{\mathrm{sing}})$.
\end{dfn}

\begin{dfn}[{\cite{BGS95}, \cite[Definition 2.1]{CMDY25}}]\label{BGS ess}
Let $X$ be a normal toric variety.
We say that a toric divisor over $X$ is \emph{BGS essential} if it appears on every toric resolution of $X$ as a prime divisor.
\end{dfn}

\begin{rem}\label{ess toric}
Let $X$ be a $\Q$-factorial affine toric variety.
By \cite[Proposition A.8]{CMDY25}, every essential divisor over $X$ is a (toric) BGS essential divisor. More precisely, a toric divisor over $X$ is essential if and only if it is exceptional and BGS essential.
\end{rem}

\begin{lem}[{cf.\ \cite[Remark 5.13]{CMKY26}}]\label{conifold ess}
Let $X$ be a $\Q$-factorial affine toric variety, and let $Y\to X$ be a proper birational toric morphism from a normal toric variety $Y$ with at worst conifold singularities.
Then every essential divisor over $X$ appears on $Y$ as a prime divisor.
\end{lem}

\begin{proof}
Each singular cone in the fan of $Y$ is a cone over a unimodular square. We subdivide each of them along one of its diagonals.
This gives a toric resolution $Y'\to Y$, which is small since no new rays are added.
Then $Y'\to X$ is a toric resolution of $X$.
Let $E$ be an essential divisor over $X$.
By Remark \ref{ess toric}, $E$ is BGS essential, so $E$ appears on $Y'$ as a prime divisor.
Since $Y'\to Y$ is small, $E$ also appears on $Y$ as a prime divisor.
\end{proof}

\section{\texorpdfstring{The case of type $\frac{1}{r}(1,r-1,a)$ with $\gcd(a,r)>1$}{The case of type 1/r(1,r-1,a) with gcd(a,r)>1}}\label{sec:type3}

In this section, let $G$ be of type $\frac{1}{r}(1,r-1,a)$ with $\gcd(a,r)>1$ and $0<a<r$, that is, of type (3) in Proposition \ref{canonical classification}.
Exchanging the first two coordinates transforms $G$ into the group of type $\frac{1}{r}(1,r-1,r-a)$, which satisfies the same assumptions. Hence every statement below has a symmetric counterpart under $x\leftrightarrow y$ (and $a\leftrightarrow r-a$).
Following \cite{Ked11}, we classify the $G$-graphs according to their number of valleys (Sections \ref{subsec:valleys3} and \ref{subsec:classification3}), compute $k[S(\Gamma)]$ in each case (Section \ref{subsec:charts3}), and deduce Theorem \ref{thm:intro conifold} for $G$ (Corollary \ref{Hilb conifold}).

\subsection{\texorpdfstring{Valleys (See \cite[Section 2]{Ked11})}{Valleys}}\label{subsec:valleys3}

We first observe that every $G$-graph lies in the union of the $xz$-plane and the $yz$-plane, and introduce the notion of valleys.

\begin{rem}
Since $\wt(xy)=\wt(1)$, we have $xy\notin\Gamma$ for any $G$-graph $\Gamma$.
Hence, every monomial in $\Gamma$ is of the form $x^mz^n$ or $y^mz^n$ for some $m,n\geq0$.
\end{rem}

\begin{dfn}\label{index}
For any $G$-graph $\Gamma$, define $\lambda(\Gamma),\mu(\Gamma),\nu(\Gamma)$ to be the unique non-negative integers such that
\begin{align*}
&x^{\lambda(\Gamma)}\in\Gamma,\ x^{\lambda(\Gamma)+1}\notin\Gamma, \\
&y^{\mu(\Gamma)}\in\Gamma,\ y^{\mu(\Gamma)+1}\notin\Gamma, \\
&z^{\nu(\Gamma)}\in\Gamma,\ z^{\nu(\Gamma)+1}\notin\Gamma.
\end{align*}
When no confusion can arise, we simply write $\lambda,\mu,\nu$.
\end{dfn}

\begin{dfn}[See {\cite[Definition 2.4]{Ked11}}]\label{xy-valley}
A monomial $x^mz^n$ (resp.\ $y^mz^n$) for $m,n\geq0$ is called an $x$-valley (resp.\ $y$-valley) of $\Gamma$ if
\begin{align*}
x^mz^n,x^{m+1}z^n,x^mz^{n+1}\in\Gamma\text{ but }x^{m+1}z^{n+1}\notin\Gamma, \\
(\text{resp.\ }y^mz^n,y^{m+1}z^n,y^mz^{n+1}\in\Gamma\text{ but }y^{m+1}z^{n+1}\notin\Gamma).
\end{align*}
We call an $x$-valley or a $y$-valley a valley for brevity.
\end{dfn}

\subsection{\texorpdfstring{Classification of $G$-graphs}{Classification of G-graphs}}\label{subsec:classification3}

In this subsection, we show that a $G$-graph has at most two valleys and determine its shape in each case.
We begin with a lemma on the weight map.

\begin{lem}\label{fund inv}
Let $\Gamma$ be a $G$-graph.
Let $v\in M_0^+\setminus\Gamma$ be such that $xy\nmid v$.
If $x^{-1}\cdot v\in\Gamma$ (resp.\ $y^{-1}\cdot v\in\Gamma$), then $\wt_\Gamma(v)$ is not divisible by $x$ (resp.\ $y$), and $y\cdot\wt_\Gamma(v)\notin\Gamma$ (resp.\ $x\cdot\wt_\Gamma(v)\notin\Gamma$).
\end{lem}

\begin{proof}
By symmetry, we may assume that $x^{-1}\cdot v\in\Gamma$.
By Lemma \ref{not divisible}, $\wt_\Gamma(v)$ is not divisible by $x$.
Since $\wt(x^{-1})=\wt(y)$, we have $\wt(x^{-1}\cdot v)=\wt(y\cdot\wt_\Gamma(v))$.
If $y\cdot\wt_\Gamma(v)\in\Gamma$, then $x^{-1}\cdot v=y\cdot\wt_\Gamma(v)$, that is, $v=xy\wt_\Gamma(v)$, which contradicts $xy\nmid v$.
\end{proof}

\begin{lem}\label{valleys}
If $x^mz^n$ is an $x$-valley of a $G$-graph $\Gamma$, then
\[
\wt_\Gamma(x^{m+1}z^{n+1})=y^{\mu(\Gamma)}.
\]
If $y^mz^n$ is a $y$-valley of $\Gamma$, then
\[
\wt_\Gamma(y^{m+1}z^{n+1})=x^{\lambda(\Gamma)}.
\]
Consequently, $\Gamma$ has at most one $x$-valley and at most one $y$-valley.
\end{lem}

\begin{proof}
Suppose that $x^mz^n$ is an $x$-valley of $\Gamma$, and set $v\coloneqq x^{m+1}z^{n+1}\notin\Gamma$.
Since $v/x,v/z\in\Gamma$ and $xy\nmid v$, Lemmas \ref{not divisible} and \ref{fund inv} show that $\wt_\Gamma(v)$ is divisible by neither $x$ nor $z$, and that $y\cdot\wt_\Gamma(v)\notin\Gamma$.
Hence $\wt_\Gamma(v)=y^j$ with $y^{j+1}\notin\Gamma$, that is, $\wt_\Gamma(v)=y^{\mu(\Gamma)}$.

Suppose that $\Gamma$ has two distinct $x$-valleys $x^mz^n$ and $x^{m'}z^{n'}$, say with $m<m'$. Then $n>n'$, as $\Gamma$ is closed under taking divisors.
By the above, $\wt(x^{m+1}z^{n+1})=\wt(x^{m'+1}z^{n'+1})$, so $\wt(z^{n-n'})=\wt(x^{m'-m})$, and hence
\[
\wt(x^mz^{n+1})=\wt(x^mz^{n'+1}z^{n-n'})=\wt(x^{m'}z^{n'+1}).
\]
Since $x^mz^{n+1}$ and $x^{m'}z^{n'+1}$ are distinct elements of $\Gamma$, this is a contradiction.
Therefore, $\Gamma$ has at most one $x$-valley.
The assertions for $y$-valleys follow by symmetry.
\end{proof}

From now on, for a fixed $G$-graph $\Gamma$, we denote its $x$-valley and its $y$-valley, if they exist, by $x^{\lambda_x}z^{\nu_x}$ and $y^{\mu_y}z^{\nu_y}$, respectively.

\begin{cor}\label{valley position}
If $x^{\lambda_x}z^{\nu_x}$ is an $x$-valley of a $G$-graph $\Gamma$, then
\[
\wt_\Gamma(y^{\mu(\Gamma)+1})=x^{\lambda_x}z^{\nu_x+1}.
\]
If $y^{\mu_y}z^{\nu_y}$ is a $y$-valley of $\Gamma$, then
\[
\wt_\Gamma(x^{\lambda(\Gamma)+1})=y^{\mu_y}z^{\nu_y+1}.
\]
\end{cor}

\begin{proof}
Suppose that $x^{\lambda_x}z^{\nu_x}$ is an $x$-valley.
By Lemma \ref{valleys}, $\wt_\Gamma(x^{\lambda_x+1}z^{\nu_x+1})=y^{\mu(\Gamma)}$.
Since
\[
\wt(y^{\mu(\Gamma)+1})=\wt(x^{\lambda_x+1}yz^{\nu_x+1})=\wt(x^{\lambda_x}z^{\nu_x+1}),
\]
and $x^{\lambda_x}z^{\nu_x+1}\in\Gamma$, we obtain $\wt_\Gamma(y^{\mu(\Gamma)+1})=x^{\lambda_x}z^{\nu_x+1}$.
The case of a $y$-valley follows by symmetry.
\end{proof}

\begin{rem}\label{shape}
By Lemma \ref{valleys}, $\Gamma$ has zero, one or two valleys.
In each case, the shape of $\Gamma$ follows from Definitions \ref{index} and \ref{xy-valley} and the fact that $\Gamma$ is closed under taking divisors.
The formulas \eqref{ijk formula1}--\eqref{ijk formula4} then follow by counting $|\Gamma|=r$ in Figures \ref{fig:zerovalley}--\ref{fig:twovalleys}.
In these figures, the shaded boxes show the monomials that span $\Gamma$.
\begin{figure}
\centering
\begin{tikzpicture}[scale=0.88]
\coordinate (A) at (-3,0) node at (A) [below]{};
\coordinate (B) at (9,0) node at (B) [below]{};
\coordinate (C) at (9,5) node at (C) [above]{};
\coordinate (D) at (-3,5) node at (D) [above]{};

\draw (A)--(D)--(C)--(B)--(A);

\coordinate  (zk) at (2.5,5) node at (zk) [above]{$z^{\nu+1}$};
\coordinate  (xi) at (9.5,0) node at (xi) [above]{$x^{\lambda+1}$};
\coordinate  (yj) at (-3.5,0) node at (yj) [above]{$y^{\mu+1}$};
\coordinate  (O) at (2.4,0)  node at (O) [above]{$1$};
\coordinate  (z) at (2.4,0.7)  node at (z) [above]{$z$};
\coordinate  (x) at (3.2,0)  node at (x) [above]{$x$};
\coordinate  (y) at (1.6,0)  node at (y) [above]{$y$};

\coordinate (1) at (2,5);
\coordinate (2) at (2,0);
\coordinate (3) at (2.8,5);
\coordinate (4) at (2.8,0);
\coordinate (5) at (2,0.6);
\coordinate (6) at (2.8,0.6);
\coordinate (7) at (2,1.2);
\coordinate (8) at (2.8,1.2);

\coordinate (9) at (-3,0.6);
\coordinate (10) at (9,0.6);
\coordinate (11) at (1.2,0);
\coordinate (12) at (3.6,0);
\coordinate (13) at (1.2,0.6);
\coordinate (14) at (3.6,0.6);

\draw[line width=1pt, dotted]  (1)--(2);
\draw[line width=1pt, dotted]  (3)--(4);
\draw[line width=1pt, dotted]  (7)--(8);
\draw[line width=1pt, dotted]  (9)--(10);
\draw[line width=1pt, dotted]  (13)--(11);
\draw[line width=1pt, dotted]  (12)--(14);

\coordinate (15) at (-3,2.4);
\coordinate (16) at (-2.4,2.4);
\coordinate (17) at (-2.4,3);

\coordinate (18) at (-3,4.4);
\coordinate (19) at (-2.4,4.4);
\coordinate (20) at (-2.4,5);

\draw[line width=1pt]  (D)--(18)--(19)--(20)--(D);
 \fill[lightgray]   (D)--(18)--(19)--(20)--(D);

 \coordinate (21) at (9,4.4);
\coordinate (22) at (8.4,4.4);
\coordinate (23) at (8.4,5);

\draw[line width=1pt]  (C)--(21)--(22)--(23)--(C);
 \fill[lightgray]   (C)--(21)--(22)--(23)--(C);

 \coordinate (24) at (9,1.4);
\coordinate (25) at (8.4,1.4);
\coordinate (26) at (8.4,2);

\end{tikzpicture}
\caption{Zero valleys}
\label{fig:zerovalley}
\end{figure}
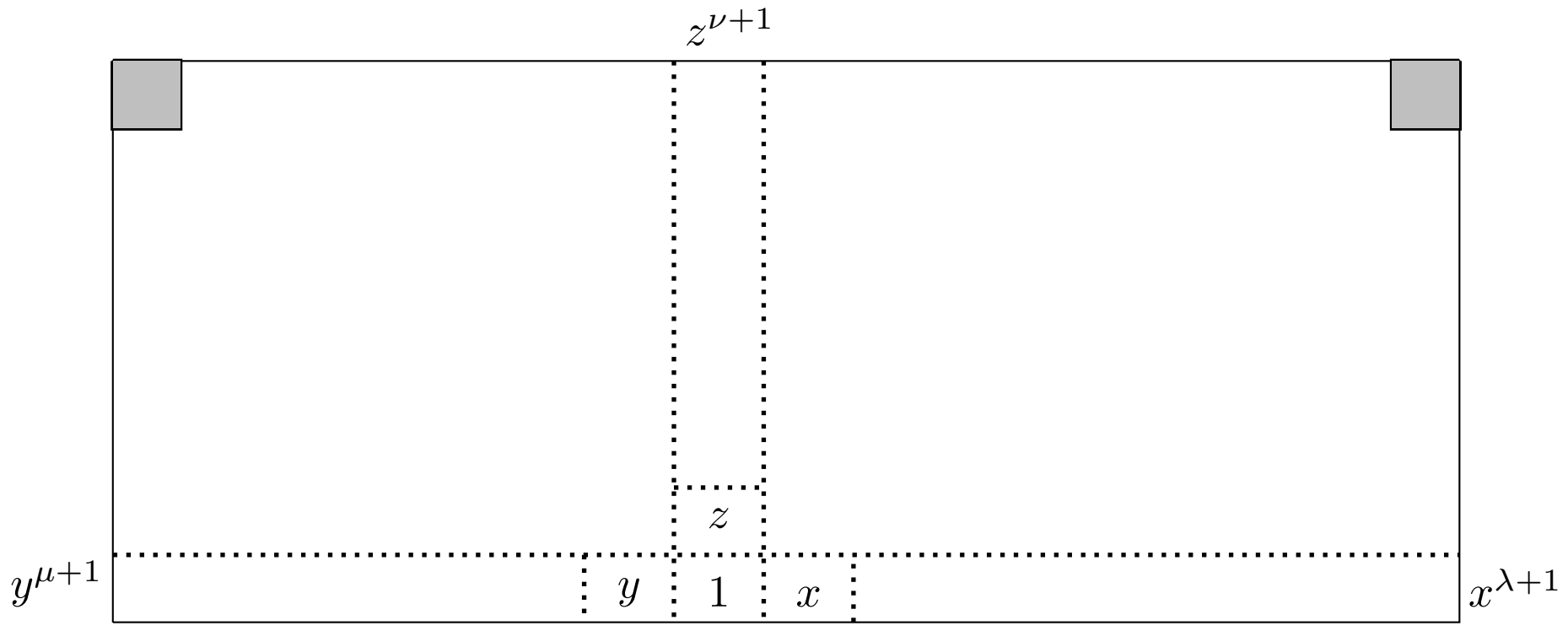

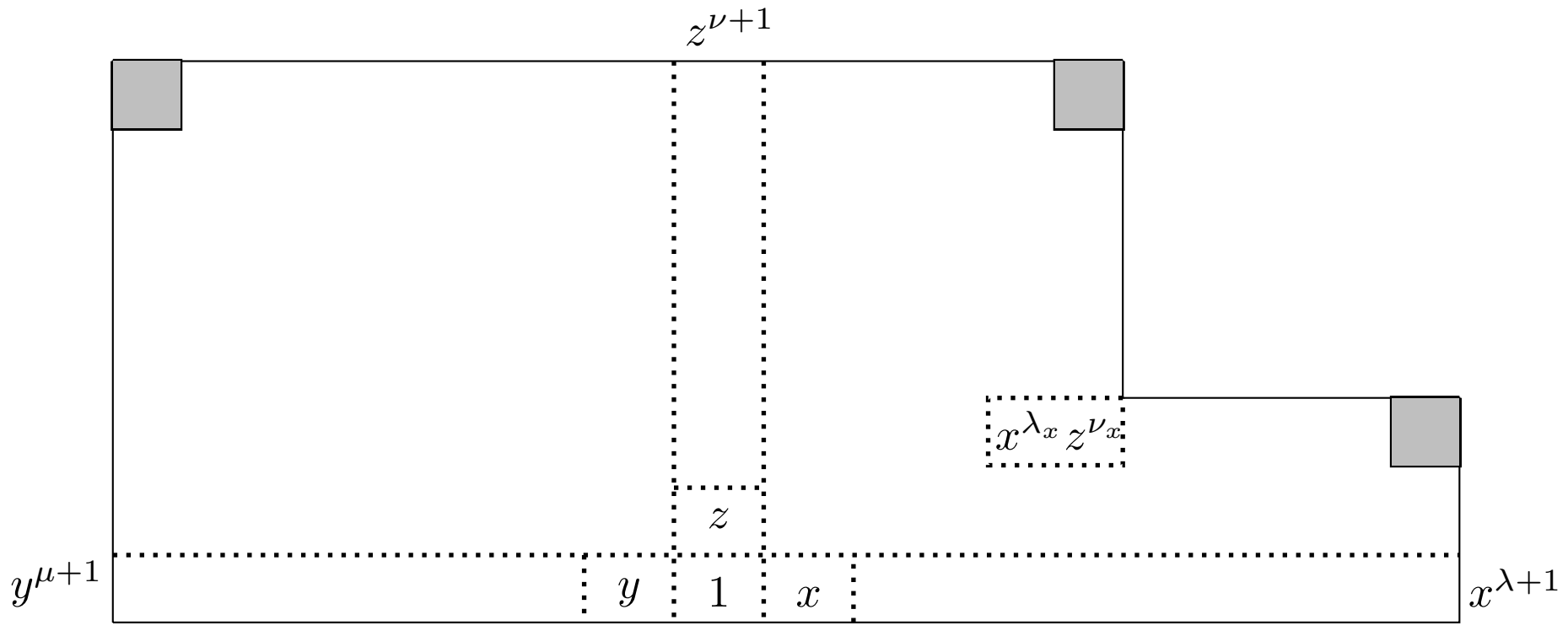
\begin{figure}
\centering
\begin{tikzpicture}[scale=0.88]
\coordinate (A) at (-3,0) node at (A) [below]{};
\coordinate (B) at (9,0) node at (B) [below]{};
\coordinate (C) at (6,5) node at (C) [above]{};
\coordinate (D) at (-3,5) node at (D) [above]{};
\coordinate  (F) at (6,2);
\coordinate  (H) at (9,2);

\draw (A)--(D)--(C)--(F)--(H)--(B)--(A);

\coordinate  (zk) at (2.5,5) node at (zk) [above]{$z^{\nu+1}$};
\coordinate  (xi) at (9.5,0) node at (xi) [above]{$x^{\lambda+1}$};
\coordinate  (yj) at (-3.5,0) node at (yj) [above]{$y^{\mu+1}$};
\coordinate  (O) at (2.4,0)  node at (O) [above]{$1$};
\coordinate  (z) at (2.4,0.7)  node at (z) [above]{$z$};
\coordinate  (x) at (3.2,0)  node at (x) [above]{$x$};
\coordinate  (y) at (1.6,0)  node at (y) [above]{$y$};

\coordinate (1) at (2,5);
\coordinate (2) at (2,0);
\coordinate (3) at (2.8,5);
\coordinate (4) at (2.8,0);
\coordinate (5) at (2,0.6);
\coordinate (6) at (2.8,0.6);
\coordinate (7) at (2,1.2);
\coordinate (8) at (2.8,1.2);

\coordinate (9) at (-3,0.6);
\coordinate (10) at (9,0.6);
\coordinate (11) at (1.2,0);
\coordinate (12) at (3.6,0);
\coordinate (13) at (1.2,0.6);
\coordinate (14) at (3.6,0.6);

\draw[line width=1pt, dotted]  (1)--(2);
\draw[line width=1pt, dotted]  (3)--(4);
\draw[line width=1pt, dotted]  (7)--(8);
\draw[line width=1pt, dotted]  (9)--(10);
\draw[line width=1pt, dotted]  (13)--(11);
\draw[line width=1pt, dotted]  (12)--(14);

\coordinate (d) at (4.8,2);
\coordinate (e) at (4.8,1.4);
\coordinate (f) at (6,1.4);
\coordinate (xv) at (5.45,1.4) node at (xv)[above]{$x^{\lambda_x}z^{\nu_x}$};

\draw[line width=1pt, dotted]  (F)--(d)--(e)--(f)--(F);

\coordinate (18) at (-3,4.4);
\coordinate (19) at (-2.4,4.4);
\coordinate (20) at (-2.4,5);

\draw[line width=1pt]  (D)--(18)--(19)--(20)--(D);
 \fill[lightgray]   (D)--(18)--(19)--(20)--(D);

 \coordinate (21) at (6,4.4);
\coordinate (22) at (5.4,4.4);
\coordinate (23) at (5.4,5);

\draw[line width=1pt]  (C)--(21)--(22)--(23)--(C);
 \fill[lightgray]   (C)--(21)--(22)--(23)--(C);

 \coordinate (24) at (9,1.4);
\coordinate (25) at (8.4,1.4);
\coordinate (26) at (8.4,2);

\draw[line width=1pt]  (H)--(24)--(25)--(26)--(H);
 \fill[lightgray] (H)--(24)--(25)--(26)--(H);

\end{tikzpicture}
\caption{One valley ($x$-valley)}
\label{fig:onevalley}
\end{figure}

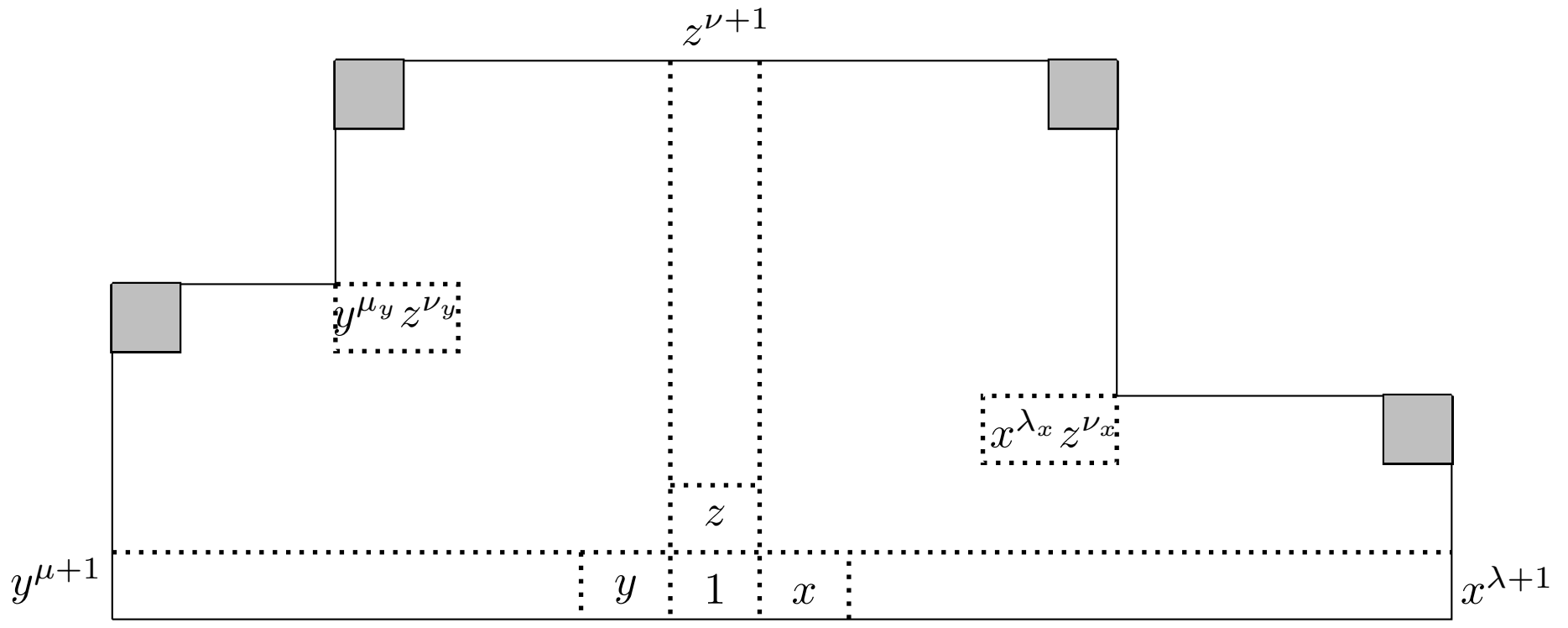
\begin{figure}
\centering
\begin{tikzpicture}[scale=0.88]
\coordinate (A) at (-3,0) node at (A) [below]{};
\coordinate (B) at (9,0) node at (B) [below]{};
\coordinate (C) at (6,5) node at (C) [above]{};
\coordinate (D) at (-1,5) node at (D) [above]{};
\coordinate  (E) at (-1,3);
\coordinate  (F) at (6,2);
\coordinate  (G) at (-3,3);
\coordinate  (H) at (9,2);

\draw (A)--(G)--(E)--(D)--(C)--(F)--(H)--(B)--(A);

\coordinate  (zk) at (2.5,5) node at (zk) [above]{$z^{\nu+1}$};
\coordinate  (xi) at (9.5,0) node at (xi) [above]{$x^{\lambda+1}$};
\coordinate  (yj) at (-3.5,0) node at (yj) [above]{$y^{\mu+1}$};
\coordinate  (O) at (2.4,0)  node at (O) [above]{$1$};
\coordinate  (z) at (2.4,0.7)  node at (z) [above]{$z$};
\coordinate  (x) at (3.2,0)  node at (x) [above]{$x$};
\coordinate  (y) at (1.6,0)  node at (y) [above]{$y$};

\coordinate (1) at (2,5);
\coordinate (2) at (2,0);
\coordinate (3) at (2.8,5);
\coordinate (4) at (2.8,0);
\coordinate (5) at (2,0.6);
\coordinate (6) at (2.8,0.6);
\coordinate (7) at (2,1.2);
\coordinate (8) at (2.8,1.2);

\coordinate (9) at (-3,0.6);
\coordinate (10) at (9,0.6);
\coordinate (11) at (1.2,0);
\coordinate (12) at (3.6,0);
\coordinate (13) at (1.2,0.6);
\coordinate (14) at (3.6,0.6);

\draw[line width=1pt, dotted]  (1)--(2);
\draw[line width=1pt, dotted]  (3)--(4);
\draw[line width=1pt, dotted]  (7)--(8);
\draw[line width=1pt, dotted]  (9)--(10);
\draw[line width=1pt, dotted]  (13)--(11);
\draw[line width=1pt, dotted]  (12)--(14);

\coordinate (a) at (-1,2.4);
\coordinate (b) at (0.1,2.4);
\coordinate (c) at (0.1,3);
\coordinate (yv) at (-0.45,2.4) node at (yv)[above]{$y^{\mu_y}z^{\nu_y}$};

\draw[line width=1pt, dotted]  (E)--(a)--(b)--(c)--(E);

\coordinate (d) at (4.8,2);
\coordinate (e) at (4.8,1.4);
\coordinate (f) at (6,1.4);
\coordinate (xv) at (5.45,1.4) node at (xv)[above]{$x^{\lambda_x}z^{\nu_x}$};

\draw[line width=1pt, dotted]  (F)--(d)--(e)--(f)--(F);

\coordinate (15) at (-3,2.4);
\coordinate (16) at (-2.4,2.4);
\coordinate (17) at (-2.4,3);

\draw[line width=1pt]  (G)--(15)--(16)--(17)--(G);
 \fill[lightgray]   (G)--(15)--(16)--(17)--(G);

\coordinate (18) at (-1,4.4);
\coordinate (19) at (-0.4,4.4);
\coordinate (20) at (-0.4,5);

\draw[line width=1pt]  (D)--(18)--(19)--(20)--(D);
 \fill[lightgray]   (D)--(18)--(19)--(20)--(D);

 \coordinate (21) at (6,4.4);
\coordinate (22) at (5.4,4.4);
\coordinate (23) at (5.4,5);

\draw[line width=1pt]  (C)--(21)--(22)--(23)--(C);
 \fill[lightgray]   (C)--(21)--(22)--(23)--(C);

 \coordinate (24) at (9,1.4);
\coordinate (25) at (8.4,1.4);
\coordinate (26) at (8.4,2);

\draw[line width=1pt]  (H)--(24)--(25)--(26)--(H);
 \fill[lightgray] (H)--(24)--(25)--(26)--(H);

\end{tikzpicture}
\caption{Two valleys}
\label{fig:twovalleys}
\end{figure}
If a $G$-graph $\Gamma$ has zero valleys, then $\Gamma=\Span(x^{\lambda}z^\nu,y^{\mu}z^\nu)$ (see Figure \ref{fig:zerovalley}).
In this case,
\begin{equation}\label{ijk formula1}
(\lambda+\mu+1)(\nu+1)=r.
\end{equation}
If a $G$-graph $\Gamma$ has one valley, then
\begin{equation*}
\Gamma=
\begin{cases*}
\Span(x^{\lambda}z^{\nu_x},x^{\lambda_x}z^\nu,y^{\mu}z^\nu) & if $\Gamma$ has an $x$-valley $x^{\lambda_x}z^{\nu_x}$, \\
\Span(x^{\lambda}z^\nu,y^{\mu_y}z^\nu,y^{\mu}z^{\nu_y}) & if $\Gamma$ has a $y$-valley $y^{\mu_y}z^{\nu_y}$
\end{cases*}
\end{equation*}
(see Figure \ref{fig:onevalley}).
If $\Gamma=\Span(x^{\lambda}z^{\nu_x},x^{\lambda_x}z^\nu,y^{\mu}z^\nu)$, then
\begin{equation}\label{ijk formula2}
(\lambda_x+\mu+1)(\nu+1)+(\lambda-\lambda_x)(\nu_x+1)=r.
\end{equation}
If $\Gamma=\Span(x^{\lambda}z^\nu,y^{\mu_y}z^\nu,y^{\mu}z^{\nu_y})$, then
\begin{equation}\label{ijk formula3}
(\lambda+\mu_y+1)(\nu+1)+(\mu-\mu_y)(\nu_y+1)=r.
\end{equation}
If a $G$-graph $\Gamma$ has two valleys $x^{\lambda_x}z^{\nu_x}$ and $y^{\mu_y}z^{\nu_y}$, then
\[
\Gamma=\Span(x^{\lambda}z^{\nu_x},x^{\lambda_x}z^\nu,y^{\mu_y}z^\nu,y^{\mu}z^{\nu_y})
\]
(see Figure \ref{fig:twovalleys}).
In this case,
\begin{equation}\label{ijk formula4}
(\lambda_x+\mu_y+1)(\nu+1)+(\lambda-\lambda_x)(\nu_x+1)+(\mu-\mu_y)(\nu_y+1)=r.
\end{equation}
\end{rem}

We now determine the $G$-graphs, and the weights of the monomials just outside them, in each of the three cases.

\begin{lem}\label{no valley}
The only possible $G$-graphs with zero valleys are
\begin{align*}
\Gamma_s^{xyz}&=\Span(x^sz^{(r/h)-1},y^{h-s-1}z^{(r/h)-1})\text{ for }s=0,1,\ldots,h-1, \\
\Gamma_s^{xy}&=\Span(x^s,y^{r-s-1})\text{ for }s=0,1,\ldots,r-1,
\end{align*}
where $h=\gcd(a,r)$.
\end{lem}

\begin{proof}
Let $\Gamma$ be a $G$-graph with no valleys, so that $\Gamma=\Span(x^\lambda z^\nu,y^\mu z^\nu)$ by Remark \ref{shape}.
Note that $\wt(x)$ has order $r$ and $\wt(z)$ has order $r/h$. In particular, $z^{r/h}\notin\Gamma$, so $\nu\leq(r/h)-1$.
Since $x^{-1}\cdot x^{\lambda+1}\in\Gamma$, Lemma \ref{fund inv} shows that $\wt_\Gamma(x^{\lambda+1})=y^tz^{t'}$ for some $t,t'\geq0$ with $y^{t+1}z^{t'}\notin\Gamma$. As $t'\leq\nu$, we have $t=\mu$.

Suppose first that $t'=0$.
Then $\wt(x^{\lambda+1})=\wt(y^\mu)=\wt(x^{-\mu})$, so $r$ divides $\lambda+\mu+1$.
Since $\lambda+\mu+1\leq r$ by \eqref{ijk formula1}, we obtain $\lambda+\mu+1=r$ and $\nu=0$, that is, $\Gamma=\Gamma^{xy}_\lambda$.

Suppose next that $t'>0$.
Since $y^{-1}\cdot y^{\mu+1}z^{t'-1}\in\Gamma$, Lemma \ref{fund inv} shows similarly that $\wt_\Gamma(y^{\mu+1}z^{t'-1})=x^\lambda z^{s}$ for some $0\leq s\leq\nu$.
Using $\wt(xy)=\wt(1)$ and $\wt(x^{\lambda+1})=\wt(y^\mu z^{t'})$, we obtain
\[
\wt(y^{\mu}z^{t'-1})=\wt(xy^{\mu+1}z^{t'-1})=\wt(x^{\lambda+1}z^{s})=\wt(y^{\mu}z^{t'+s}),
\]
so $\wt(z^{s+1})=\wt(1)$, and hence $r/h$ divides $s+1$.
Since $s\leq\nu\leq(r/h)-1$, we obtain $s=\nu=(r/h)-1$, and \eqref{ijk formula1} gives $\lambda+\mu+1=h$, that is, $\Gamma=\Gamma^{xyz}_\lambda$.
\end{proof}

\begin{lem}\label{one valley}
Let $\Gamma$ be a $G$-graph whose only valley is an $x$-valley $x^{\lambda_x}z^{\nu_x}$.
Then
\[
\wt_\Gamma(x^{\lambda+1})=y^{\mu}z^{\nu-\nu_x},\quad\wt_\Gamma(z^{\nu+1})=x^{\lambda-\lambda_x}.
\]
The symmetric statement holds if the only valley of $\Gamma$ is a $y$-valley.
\end{lem}

\begin{proof}
We have $\Gamma=\Span(x^{\lambda}z^{\nu_x},x^{\lambda_x}z^\nu,y^{\mu}z^\nu)$ by Remark \ref{shape}. In particular, $y^jz^q\in\Gamma$ if and only if $j\leq\mu$ and $q\leq\nu$.
Recall also that $\wt(z)$ has order $r/h$, so that $\nu\leq(r/h)-1$.
For $0\leq m\leq \nu_x$, since $x^{-1}\cdot x^{\lambda+1}z^m\in\Gamma$, Lemma \ref{fund inv} shows that $\wt_\Gamma(x^{\lambda+1}z^m)=y^{t}z^{t'}$ for some $t,t'\geq0$ with $y^{t+1}z^{t'}\notin\Gamma$. As $t'\leq\nu$, we have $t=\mu$.
Thus, for any $0\leq m\leq \nu_x$, we have
\begin{equation}\label{x-valley}
\wt_\Gamma(x^{\lambda+1}z^m)=y^{\mu}z^{l_m}
\end{equation}
for some $l_m\geq0$.
Since $\wt_\Gamma(x^{\lambda_x+1}z^{\nu_x+1})=y^\mu$ (Lemma \ref{valleys}), we have
\begin{align*}
\wt(x^{\lambda}z^m)&=\wt(x^{\lambda+1}yz^m)=\wt(y^{\mu+1}z^{l_m}) \\
&=\wt(x^{\lambda_x+1}yz^{\nu_x+l_m+1})=\wt(x^{\lambda_x}z^{\nu_x+l_m+1}).
\end{align*}
If $l_m<\nu-\nu_x$ for some $m$, then $x^{\lambda}z^m$ and $x^{\lambda_x}z^{\nu_x+l_m+1}$ are distinct elements of $\Gamma$ of the same weight, a contradiction. Hence $l_m\geq \nu-\nu_x$ for all $0\leq m\leq \nu_x$.
Suppose that $l_0>\nu-\nu_x$, and set $m_0\coloneqq\nu-l_0+1$, so that $1\leq m_0\leq\nu_x$.
Since $\wt(x^{\lambda+1})=\wt(y^\mu z^{l_0})$, we have
\[
\wt(y^{\mu}z^{l_{m_0}})=\wt(x^{\lambda+1}z^{m_0})=\wt(y^{\mu}z^{l_0+m_0})=\wt(y^\mu z^{\nu+1}),
\]
so $\wt(z^{\nu+1-l_{m_0}})=\wt(1)$.
Since $1\leq\nu+1-l_{m_0}\leq\nu_x+1\leq\nu<r/h$, this contradicts the fact that $\wt(z)$ has order $r/h$.
Therefore $l_0=\nu-\nu_x$, that is, $\wt_\Gamma(x^{\lambda+1})=y^{\mu}z^{\nu-\nu_x}$.

Since $\wt(y^\mu)=\wt(x^{\lambda_x+1}z^{\nu_x+1})$ by Lemma \ref{valleys}, we have
\[
\wt(x^{\lambda+1})=\wt(y^{\mu}z^{\nu-\nu_x})=\wt(x^{\lambda_x+1}z^{\nu+1}),
\]
so $\wt(z^{\nu+1})=\wt(x^{\lambda-\lambda_x})$.
Since $x^{\lambda-\lambda_x}\in\Gamma$, we obtain $\wt_\Gamma(z^{\nu+1})=x^{\lambda-\lambda_x}$.
\end{proof}

\begin{lem}\label{two valley}
Let $\Gamma$ be a $G$-graph with two valleys $x^{\lambda_x}z^{\nu_x}$ and $y^{\mu_y}z^{\nu_y}$.
Then
\[
\nu=\nu_x+\nu_y+1
\]
and $\wt_\Gamma(z^{\nu+1})=v$, where
\begin{equation*}
v\coloneqq
\begin{cases*}
x^{(\lambda-\lambda_x)-(\mu-\mu_y)} & if $\lambda-\lambda_x\geq \mu-\mu_y$, \\
y^{(\mu-\mu_y)-(\lambda-\lambda_x)} & if $\lambda-\lambda_x<\mu-\mu_y$.
\end{cases*}
\end{equation*}
\end{lem}

\begin{proof}
By Remark \ref{shape}, $\Gamma=\Span(x^{\lambda}z^{\nu_x},x^{\lambda_x}z^\nu,y^{\mu_y}z^\nu,y^{\mu}z^{\nu_y})$, and $\wt(z)$ has order $r/h>\nu$.
By Lemma \ref{valleys}, $\wt(x^{\lambda_x+1}z^{\nu_x+1})=\wt(y^\mu)$ and $\wt(y^{\mu_y+1}z^{\nu_y+1})=\wt(x^\lambda)$.
Since $y^{-1}\cdot y^{\mu+1}z^{\nu_y}\in\Gamma$, Lemma \ref{fund inv} shows that $\wt_\Gamma(y^{\mu+1}z^{\nu_y})=x^{t}z^{t'}$ for some $t,t'\geq0$ with $x^{t+1}z^{t'}\notin\Gamma$. By the shape of $\Gamma$, either $t=\lambda$ and $t'\leq\nu_x$, or $t=\lambda_x$ and $\nu_x<t'\leq\nu$.
In the former case,
\begin{align*}
\wt(x^{\lambda-\lambda_x}y^\mu)&=\wt(x^{\lambda+1}z^{\nu_x+1})=\wt(xy^{\mu+1}z^{\nu_x+\nu_y+1-t'})=\wt(y^{\mu}z^{\nu_x+\nu_y+1-t'}) \\
&=\wt(x^{\lambda}y^{\mu-\mu_y-1}z^{\nu_x-t'}),
\end{align*}
so $\wt(x^{\lambda_x+1}z^{\nu_x-t'})=\wt(y^{\mu_y})$.
This is a contradiction, since these are distinct elements of $\Gamma$.
Hence $t=\lambda_x$, and
\[
\wt(y^{\mu}z^{\nu_y})=\wt(xy^{\mu+1}z^{\nu_y})=\wt(x^{\lambda_x+1}z^{t'})=\wt(y^{\mu}z^{t'-\nu_x-1}).
\]
Since $0\leq\nu_y,t'-\nu_x-1\leq\nu<r/h$, we obtain $t'=\nu_x+\nu_y+1$, and in particular $\nu_x+\nu_y+1\leq\nu$.

On the other hand,
\[
\wt(x^{\lambda_x}y^{\mu_y}z^{\nu_x+\nu_y+2})=\wt(x^{\lambda_x+1}z^{\nu_x+1}\cdot y^{\mu_y+1}z^{\nu_y+1})=\wt(x^{\lambda}y^\mu),
\]
so $\wt(z^{\nu_x+\nu_y+2})=\wt(x^{\lambda-\lambda_x}y^{\mu-\mu_y})=\wt(v)$.
Since $v\in\Gamma$, we have $z^{\nu_x+\nu_y+2}\notin\Gamma$, that is, $\nu=\nu_x+\nu_y+1$, and $\wt_\Gamma(z^{\nu+1})=v$.
\end{proof}

\subsection{\texorpdfstring{The coordinate rings $k[S(\Gamma)]$}{The coordinate rings k[S(Γ)]}}\label{subsec:charts3}

Using the above classification, we compute $k[S(\Gamma)]$ for every $G$-graph $\Gamma$ with $\dim\sigma(\Gamma)=3$, according to the number of valleys of $\Gamma$.

\begin{thm}\label{zero valley chart}
Let $\Gamma$ be a $G$-graph with zero valleys and $\dim\sigma(\Gamma)=3$.
Then $k[S(\Gamma)]\cong k[x,y,z]$.
\end{thm}

\begin{proof}
By Lemma \ref{no valley}, $\Gamma=\Gamma_s^{xy}$ or $\Gamma=\Gamma_s^{xyz}$ for some $s$.
By Proposition \ref{monoid generator}, $S(\Gamma)$ is generated by $s_\Gamma(x^{\lambda+1}),s_\Gamma(y^{\mu+1}),s_\Gamma(z^{\nu+1}),s_\Gamma(xy)$.
First assume that $\Gamma=\Gamma_s^{xy}$.
Then $\nu=0$, so $\lambda+\mu+1=r$ by \eqref{ijk formula1}, and hence
\begin{align*}
s_\Gamma(xy)&=xy,s_\Gamma(x^{\lambda+1})=x^{\lambda+1}y^{\lambda+1-r},s_\Gamma(y^{\mu+1})=x^{\mu+1-r}y^{\mu+1}, \\
s_\Gamma(z^{\nu+1})&=
\begin{cases*}
x^{-a}z^{\nu+1} & if $x^a\in\Gamma$, \\
y^{a-r}z^{\nu+1} & if $y^{r-a}\in\Gamma$.
\end{cases*}
\end{align*}
Since
\[
s_\Gamma(xy)=s_\Gamma(x^{\lambda+1})s_\Gamma(y^{\mu+1}),
\]
$S(\Gamma)$ is generated by $s_\Gamma(x^{\lambda+1})$, $s_\Gamma(y^{\mu+1})$, $s_\Gamma(z^{\nu+1})$, and Lemma \ref{basis} applies.

Next assume that $\Gamma=\Gamma_s^{xyz}$.
Then $\nu=(r/h)-1$, so by \eqref{ijk formula1}, we have
\begin{equation}\label{ijg formula}
\lambda+\mu+1=\frac{r}{\nu+1}=h.
\end{equation}
Since $h=\gcd(a,r)$, there is an integer $0<m<r/h$ with $am\equiv h\pmod{r}$.
By \eqref{ijg formula},
\begin{align*}
\wt(x^{\lambda+1})&=\wt(x^{h-\mu})=\wt(y^{\mu}z^m), \\
\wt(y^{\mu+1})&=\wt(y^{h-\lambda})=\wt(x^{\lambda}z^{(r/h)-m}).
\end{align*}
Since $y^{\mu}z^m,x^{\lambda}z^{(r/h)-m}\in\Gamma$, we have
\begin{align*}
s_\Gamma(x^{\lambda+1})&=x^{\lambda+1}y^{-\mu}z^{-m},s_\Gamma(y^{\mu+1})=x^{-\lambda}y^{\mu+1}z^{-(r/h)+m}, \\
s_\Gamma(z^{\nu+1})&=z^{\nu+1},s_\Gamma(xy)=xy.
\end{align*}
Since
\[
s_\Gamma(xy)=s_\Gamma(x^{\lambda+1})s_\Gamma(y^{\mu+1})s_\Gamma(z^{\nu+1}),
\]
$S(\Gamma)$ is generated by $s_\Gamma(x^{\lambda+1})$, $s_\Gamma(y^{\mu+1})$, $s_\Gamma(z^{\nu+1})$, and Lemma \ref{basis} applies.
\end{proof}

\begin{thm}\label{one valley chart}
Let $\Gamma$ be a $G$-graph with one valley and $\dim\sigma(\Gamma)=3$.
Then $k[S(\Gamma)]\cong k[x,y,z]$.
\end{thm}

\begin{proof}
By symmetry, we may assume that the valley is an $x$-valley $x^{\lambda_x}z^{\nu_x}$.
By Proposition \ref{monoid generator}, $S(\Gamma)$ is generated by
\[
s_\Gamma(x^{\lambda+1}),\ s_\Gamma(y^{\mu+1}),\ s_\Gamma(z^{\nu+1}),\ s_\Gamma(xy),\ s_\Gamma(x^{\lambda_x+1}z^{\nu_x+1}).
\]
By Corollary \ref{valley position} and Lemma \ref{one valley}, we have
\begin{align*}
s_\Gamma(x^{\lambda+1})&=x^{\lambda+1}y^{-\mu}z^{-\nu+\nu_x},s_\Gamma(y^{\mu+1})=x^{-\lambda_x}y^{\mu+1}z^{-\nu_x-1}, \\
s_\Gamma(z^{\nu+1})&=x^{-\lambda+\lambda_x}z^{\nu+1},s_\Gamma(x^{\lambda_x+1}z^{\nu_x+1})=x^{\lambda_x+1}y^{-\mu}z^{\nu_x+1}, \\
s_\Gamma(xy)&=xy.
\end{align*}
Since
\begin{align*}
s_\Gamma(xy)&=s_\Gamma(x^{\lambda+1})s_\Gamma(y^{\mu+1})s_\Gamma(z^{\nu+1}), \\
s_\Gamma(x^{\lambda_x+1}z^{\nu_x+1})&=s_\Gamma(x^{\lambda+1})s_\Gamma(z^{\nu+1}),
\end{align*}
$S(\Gamma)$ is generated by $s_\Gamma(x^{\lambda+1})$, $s_\Gamma(y^{\mu+1})$, $s_\Gamma(z^{\nu+1})$, and Lemma \ref{basis} applies.
\end{proof}

\begin{thm}\label{two valley chart}
Let $\Gamma$ be a $G$-graph with two valleys and $\dim\sigma(\Gamma)=3$.
Then $k[S(\Gamma)]\cong k[x,y,z,w]/(xy-zw)$.
\end{thm}

\begin{proof}
By Proposition \ref{monoid generator}, $S(\Gamma)$ is generated by $s_\Gamma(x^{\lambda+1})$, $s_\Gamma(y^{\mu+1})$, $s_\Gamma(z^{\nu+1})$, $s_\Gamma(xy)$, $s_\Gamma(x^{\lambda_x+1}z^{\nu_x+1})$, $s_\Gamma(y^{\mu_y+1}z^{\nu_y+1})$.
By Corollary \ref{valley position} and Lemma \ref{two valley}, we have $\nu=\nu_x+\nu_y+1$,
\begin{align*}
s_\Gamma(x^{\lambda+1})&=x^{\lambda+1}y^{-\mu_y}z^{-\nu_y-1},s_\Gamma(y^{\mu+1})=x^{-\lambda_x}y^{\mu+1}z^{-\nu_x-1}, \\
s_\Gamma(z^{\nu+1})&=
\begin{cases*}
x^{(\mu-\mu_y)-(\lambda-\lambda_x)}z^{\nu+1} & if $\lambda-\lambda_x\geq \mu-\mu_y$, \\
y^{(\lambda-\lambda_x)-(\mu-\mu_y)}z^{\nu+1} & if $\lambda-\lambda_x<\mu-\mu_y$,
\end{cases*}\\
s_\Gamma(xy)&=xy,s_\Gamma(x^{\lambda_x+1}z^{\nu_x+1})=x^{\lambda_x+1}y^{-\mu}z^{\nu_x+1}, \\
s_\Gamma(y^{\mu_y+1}z^{\nu_y+1})&=x^{-\lambda}y^{\mu_y+1}z^{\nu_y+1}.
\end{align*}
Since
\begin{align*}
s_\Gamma(xy)&=s_\Gamma(y^{\mu+1})s_\Gamma(x^{\lambda_x+1}z^{\nu_x+1})=s_\Gamma(x^{\lambda+1})s_\Gamma(y^{\mu_y+1}z^{\nu_y+1}), \\
s_\Gamma(z^{\nu+1})&=
\begin{cases*}
s_\Gamma(x^{\lambda_x+1}z^{\nu_x+1})s_\Gamma(y^{\mu_y+1}z^{\nu_y+1})s_\Gamma(xy)^{\mu-\mu_y-1} & if $\lambda-\lambda_x\geq \mu-\mu_y$, \\
s_\Gamma(x^{\lambda_x+1}z^{\nu_x+1})s_\Gamma(y^{\mu_y+1}z^{\nu_y+1})s_\Gamma(xy)^{\lambda-\lambda_x-1} & if $\lambda-\lambda_x<\mu-\mu_y$,
\end{cases*}
\end{align*}
$S(\Gamma)$ is generated by $s_\Gamma(x^{\lambda+1})$, $s_\Gamma(y^{\mu+1})$, $s_\Gamma(x^{\lambda_x+1}z^{\nu_x+1})$, $s_\Gamma(y^{\mu_y+1}z^{\nu_y+1})$, which satisfy $s_\Gamma(y^{\mu+1})s_\Gamma(x^{\lambda_x+1}z^{\nu_x+1})=s_\Gamma(x^{\lambda+1})s_\Gamma(y^{\mu_y+1}z^{\nu_y+1})$, and Lemma \ref{basis} applies.
\end{proof}

Combining Theorems \ref{zero valley chart}, \ref{one valley chart} and \ref{two valley chart}, we obtain Theorem \ref{thm:intro conifold} for groups of type (3).

\begin{cor}\label{Hilb conifold}
Let $G$ be of type $\frac{1}{r}(1,r-1,a)$ with $\gcd(a,r)>1$ and $0<a<r$, and let $\Gamma$ be a $G$-graph with $\dim\sigma(\Gamma)=3$.
Then $V(\Gamma)=U_{\sigma(\Gamma)}$.
Moreover, $\Hilb^G(\A_k^3)$ has at worst conifold singularities.
In particular, $\Hilb^G(\A_k^3)$ is normal.
\end{cor}

\begin{proof}
By Lemma \ref{valleys}, every $G$-graph falls under one of Theorems \ref{zero valley chart}, \ref{one valley chart} and \ref{two valley chart}, so $k[S(\Gamma)]$ is isomorphic to $k[x,y,z]$ or $k[x,y,z,w]/(xy-zw)$.
Both are normal, so $V(\Gamma)=U_{\sigma(\Gamma)}$.
The remaining assertions follow from Theorem \ref{Nakamura fan}.
\end{proof}

\section{\texorpdfstring{The case of type $\frac{1}{4t}(1,2t+1,4t-2)$ with $t\geq2$}{The case of type 1/4t(1,2t+1,4t-2) with t≥2}}\label{sec:type4}

In this section, let $G$ be of type $\frac{1}{4t}(1,2t+1,4t-2)$ with $t\geq2$, that is, of type (4) in Proposition \ref{canonical classification}.
Let $g\coloneqq\diag(\varepsilon,\varepsilon^{2t+1},\varepsilon^{4t-2})$ be the generator of $G$.
Since $(2t+1)^2\equiv1$ and $(4t-2)(2t+1)\equiv4t-2\pmod{4t}$, the element $g^{2t+1}=\diag(\varepsilon^{2t+1},\varepsilon,\varepsilon^{4t-2})$ also generates $G$.
Hence $G$ is preserved by exchanging the first two coordinates, and every statement below is symmetric under $x\leftrightarrow y$.
Unlike in Section \ref{sec:type3}, a $G$-graph may contain $xy$, but only when it does not contain $z^{t-1}$ (Remark \ref{zl-1 in Gamma}).
We therefore treat two cases separately.
The $G$-graphs containing $z^{t-1}$ can be handled by valleys as in Section \ref{sec:type3}, and we treat them in Section \ref{subsec:z in}.
The other $G$-graphs are treated in Section \ref{subsec:z notin}.
In each case, we classify the $G$-graphs and compute $k[S(\Gamma)]$.
We then deduce Theorem \ref{thm:intro conifold} for $G$ in Corollary \ref{Hilb conifold section4}.

\subsection{A few remarks}

We first collect relations among the weights of monomials, which follow directly from the definition of $G$.

\begin{rem}\label{wt relations}
We will use the following relations without further mention:
\begin{gather*}
\wt(z^{2t})=\wt(x^2z)=\wt(y^2z)=\wt(1),\\
\wt(xz^t)=\wt(x^{2t+1})=\wt(y),\quad\wt(yz^t)=\wt(x).
\end{gather*}
For $0\leq m\leq 2t$,
\[
\wt(z^{2t-m})=\wt(x^{2m})=\wt(y^{2m}).
\]
For $0\leq n\leq t$,
\[
\wt(yz^{t-n})=\wt(x^{2n+1}),\quad\wt(xz^{t-n})=\wt(y^{2n+1}).
\]
For $0\leq m\leq t-1$,
\[
\wt(z^{t-1-m})=\wt(xy^{2m+1})=\wt(x^{2m+1}y).
\]
\end{rem}

\begin{rem}\label{zl-1 in Gamma}
We have $\wt(xy)=\wt(z^{t-1})$.
Thus, if a $G$-graph $\Gamma$ contains $z^{t-1}$, then $xy\notin\Gamma$, and hence every monomial in $\Gamma$ is of the form $x^mz^n$ or $y^mz^n$ for some $m,n\geq0$.
\end{rem}

\subsection{\texorpdfstring{Classification of $G$-graphs $\Gamma$ with $z^{t-1}\in\Gamma$}{Classification of G-graphs Γ with z\string^(t-1) ∈ Γ}}\label{subsec:z in}

We first consider $G$-graphs with $z^{t-1}\in\Gamma$.
By Remark \ref{zl-1 in Gamma}, every monomial in such a $\Gamma$ is of the form $x^mz^n$ or $y^mz^n$ for some $m,n\geq0$, as in Section \ref{sec:type3}. Hence we can use the same notation as in Definitions \ref{index} and \ref{xy-valley}.

\begin{lem}\label{z2l-1}
Let $\Gamma$ be a $G$-graph with $z^{t-1}\in\Gamma$.
If $x^2,y^2\notin\Gamma$, then $\Gamma$ is one of the following three $G$-graphs:
\[
\Gamma_1\coloneqq\Span(xz^{2t-1}),\Gamma_2\coloneqq\Span(yz^{2t-1}),\Gamma_3\coloneqq\Span(xz^{t-1},z^{2t-1},yz^{t-1}).
\]
\end{lem}

\begin{proof}
Suppose that $x^2,y^2\notin\Gamma$.
By Remark \ref{zl-1 in Gamma}, we have $xy\notin\Gamma$, and since $\wt(z^{2t})=\wt(1)$, we have $z^{2t}\notin\Gamma$.
Hence
\[
\Gamma\subset\{z^c,xz^c,yz^c\mid0\leq c\leq2t-1\}.
\]
If $y\notin\Gamma$, then $\Gamma\subset\{z^c,xz^c\mid0\leq c\leq2t-1\}$. Since the right-hand side has exactly $4t=|\Gamma|$ elements, we obtain $\Gamma=\Span(xz^{2t-1})=\Gamma_1$.
Similarly, if $x\notin\Gamma$, then $\Gamma=\Gamma_2$.
Suppose that $x,y\in\Gamma$.
Then $\Gamma=\Span(xz^a,yz^b,z^c)$ for some $0\leq a,b\leq c\leq2t-1$.
Since $\wt(xz^t)=\wt(y)$ and $\wt(yz^t)=\wt(x)$, we have $xz^t,yz^t\notin\Gamma$, that is, $a,b\leq t-1$.
Hence
\[
4t=(a+1)+(b+1)+(c+1)\leq t+t+2t=4t,
\]
so $a=b=t-1$ and $c=2t-1$, that is, $\Gamma=\Gamma_3$.
\end{proof}

Next we consider the case where $x^2\in\Gamma$ or $y^2\in\Gamma$.
Since $\wt(x^2)=\wt(y^2)$, if $x^2\in\Gamma$, then $y^2\notin\Gamma$.

\begin{lem}\label{xzyz}
Let $\Gamma$ be a $G$-graph with $z^{t-1}\in\Gamma$.
If $x^2\in\Gamma$ and $y\notin\Gamma$, then
\[
\Gamma=\Gamma^{xz}_s\coloneqq\Span(x^{2s-1},xz^{2t-s})
\]
for some $2\leq s\leq t+1$.
If $y^2\in\Gamma$ and $x\notin\Gamma$, then
\[
\Gamma=\Gamma^{yz}_s\coloneqq\Span(y^{2s-1},yz^{2t-s})
\]
for some $2\leq s\leq t+1$.
\end{lem}

\begin{proof}
By symmetry, we may assume that $x^2\in\Gamma$ and $y\notin\Gamma$.
Since $y\notin\Gamma$, every monomial in $\Gamma$ is of the form $x^iz^c$, and since $\wt(x^2z)=\wt(1)$, we have $x^2z\notin\Gamma$.
Hence we can write $\Gamma=\Span(x^a,xz^b,z^c)$ with $a\geq2$ and $0\leq b\leq c$, where $c\geq t-1$ as $z^{t-1}\in\Gamma$.
Since $\wt(z^{2t-1})=\wt(x^2)$ and $x^2\in\Gamma$, we have $z^{2t-1}\notin\Gamma$, so $c\leq2t-2$.
Write $c=2t-s$ with $2\leq s\leq t+1$.
Since $\wt(x^{2s})=\wt(z^{2t-s})$ and $z^{2t-s}\in\Gamma$, we have $x^{2s}\notin\Gamma$, so $a\leq2s-1$. Moreover, $b\leq c=2t-s$.
Counting the monomials $x^i$ ($0\leq i\leq a$), $z^q$ ($1\leq q\leq c$) and $xz^q$ ($1\leq q\leq b$), we obtain
\[
4t=|\Gamma|=(a+1)+c+b\leq2s+(2t-s)+(2t-s)=4t,
\]
so $a=2s-1$ and $b=c=2t-s$, that is, $\Gamma=\Gamma^{xz}_s$.
\end{proof}

\begin{lem}\label{xxyyyx}
Let $\Gamma$ be a $G$-graph with $z^{t-1}\in\Gamma$.
If $x^2,y\in\Gamma$, then every $x$-valley of $\Gamma$ is $x$ or of the form $z^{\nu_x}$, and every $y$-valley of $\Gamma$ is of the form $z^{\nu_y}$.
In particular, $\Gamma$ has at most two $x$-valleys and at most one $y$-valley, hence at most three valleys.
The symmetric statement holds if $x,y^2\in\Gamma$.
\end{lem}

\begin{proof}
By symmetry, we may assume that $x^2,y\in\Gamma$.

Since $\wt(x^2z)=\wt(1)$, we have $x^2z\notin\Gamma$, and hence $x^mz^n\notin\Gamma$ for all $m\geq2$ and $n\geq1$.
For $0\leq m\leq\lambda$, set $f(m)\coloneqq\max\{n\geq0\mid x^mz^n\in\Gamma\}$. Then $f$ is non-increasing, and the above gives $f(m)=0$ for all $2\leq m\leq\lambda$.
The $x$-valleys of $\Gamma$ are exactly the monomials $x^mz^{f(m+1)}$ with $0\leq m\leq\lambda-1$ and $f(m)>f(m+1)$. Since $f(m)>f(m+1)$ is possible only for $m=0$ or $m=1$, every $x$-valley of $\Gamma$ is $x^0z^{f(1)}$ or $x^1z^{f(2)}=x$.
In particular, $\Gamma$ has at most two $x$-valleys.

Since $\wt(x^2)=\wt(y^2)$ and $x^2\in\Gamma$, we have $y^2\notin\Gamma$. Together with $y\in\Gamma$, this gives $\mu=1$.
If $y^mz^n$ is a $y$-valley, then $y^{m+1}z^n\in\Gamma$, so $m+1\leq\mu=1$, that is, $m=0$. Hence every $y$-valley of $\Gamma$ is of the form $y^0z^n=z^n$.
Moreover, such a $y$-valley satisfies $yz^n\in\Gamma$ and $yz^{n+1}\notin\Gamma$, so $n=\max\{q\geq0\mid yz^q\in\Gamma\}$ is uniquely determined.
In particular, $\Gamma$ has at most one $y$-valley.
\end{proof}

\begin{lem}\label{x2lyy2lx}
Let $\Gamma$ be a $G$-graph as in Lemma \ref{xxyyyx}.
Then $\Gamma$ has at least two valleys.
If $\Gamma$ has exactly two valleys, then $\Gamma=\Gamma_4\coloneqq\Span(x^{2t},xz^{t-1},y)$ or $\Gamma=\Gamma_5\coloneqq\Span(y^{2t},yz^{t-1},x)$.
\end{lem}

\begin{proof}
By symmetry, we may assume that $x^2,y\in\Gamma$.
Let $f$ be as in the proof of Lemma \ref{xxyyyx}.
Then $x$ is an $x$-valley if and only if $f(1)>f(2)=0$, that is, $xz\in\Gamma$.
Also, $z^{f(1)}$ is an $x$-valley if and only if $f(1)<f(0)=\nu$.
Since $\nu\geq t-1\geq1$, at least one of these holds, so $\Gamma$ has an $x$-valley.
Suppose that $\Gamma$ has no $y$-valley, that is, $yz^\nu\in\Gamma$.
Since $\wt(yz^t)=\wt(x)$, we have $\nu=t-1$, and since $\wt(yz^{t-1})=\wt(x^3)$, we have $\lambda=2$.
Then every monomial in $\Gamma$ is $x^2$ or of the form $z^q,xz^q,yz^q$ with $0\leq q\leq t-1$, so $|\Gamma|\leq3t+1<4t$, a contradiction.
Hence $\Gamma$ has a $y$-valley, and in particular at least two valleys.

Now suppose that $\Gamma$ has exactly two valleys.
By the above, $\Gamma$ has exactly one $x$-valley and one $y$-valley $z^{\nu_y}$ with $\nu_y<\nu$.
By Lemma \ref{xxyyyx}, the $x$-valley is either $z^{\nu_x}$ or $x$.
In the former case $x$ is not an $x$-valley, so $xz\notin\Gamma$. In the latter case $z^{\nu_x}$ is not an $x$-valley, so $xz^\nu\in\Gamma$.
Hence $\Gamma$ is one of the following:
\[
\Span(x^\lambda,yz^{\nu_y},z^\nu),\quad\Span(x^\lambda,xz^\nu,yz^{\nu_y}).
\]
Suppose that $\Gamma=\Span(x^\lambda,yz^{\nu_y},z^\nu)$.
Since $\wt(x^{2t+1})=\wt(y)$, $\wt(yz^t)=\wt(x)$ and $\wt(z^{2t-1})=\wt(x^2)$, we have $\lambda\leq 2t$, $\nu_y\leq t-1$ and $\nu\leq 2t-2$.
Write $\nu=2t-m$ with $2\leq m\leq t+1$.
Since $\wt(z^{2t-m})=\wt(x^{2m})$, we have $\lambda\leq 2m-1$. Moreover, $\nu_y\leq\nu-1=2t-m-1$.
Counting the monomials $x^i$ ($0\leq i\leq\lambda$), $z^q$ ($1\leq q\leq\nu$) and $yz^q$ ($0\leq q\leq\nu_y$), we obtain
\[
4t=(\lambda+1)+\nu+(\nu_y+1)\leq2m+(2t-m)+(2t-m)=4t.
\]
Hence $\nu_y=2t-m-1$ and $\lambda=2m-1$.
Since $\nu_y\leq t-1$ and $\lambda\leq2t$, we obtain $m=t$, so that $\nu_y=t-1$ and $\lambda=2t-1\geq3$.
However, $\wt(yz^{t-1})=\wt(x^3)$ and $yz^{t-1},x^3\in\Gamma$, a contradiction.
Hence $\Gamma=\Span(x^\lambda,xz^\nu,yz^{\nu_y})$.
Since $\wt(xz^t)=\wt(y)$, we have $\nu=t-1$.
Since $\nu_y<\nu=t-1$, we can write $\nu_y=t-m$ with $2\leq m\leq t$.
Since $\wt(yz^{t-m})=\wt(x^{2m+1})$, we have $\lambda\leq2m$.
Counting the monomials $z^q$, $xz^q$ ($0\leq q\leq\nu$), $yz^q$ ($0\leq q\leq\nu_y$) and $x^i$ ($2\leq i\leq\lambda$), we obtain
\[
4t=2(\nu+1)+(\nu_y+1)+(\lambda-1)\leq2t+(t-m+1)+(2m-1)=3t+m\leq4t.
\]
Hence $m=t$, so $\nu_y=0$ and $\lambda=2t$.
Therefore, $\Gamma=\Span(x^{2t},xz^{t-1},y)=\Gamma_4$.
\end{proof}

\begin{lem}\label{3 valleys}
Let $\Gamma$ be a $G$-graph with $z^{t-1},x^2,y\in\Gamma$, and assume that $\Gamma$ has three valleys. By Lemma \ref{xxyyyx}, these are $x$, $z^{\nu_x}$ and $z^{\nu_y}$ for some $\nu_x,\nu_y\geq0$.
Then $\nu_x=t-1$, and hence $\wt_\Gamma(xz^{\nu_x+1})=y$.
Moreover, the following hold.
\begin{enumerate}
\item If $\nu-\nu_y<t$, then $\nu=\nu_x+\nu_y$, $\lambda=2(t-\nu_y)$ and
\[
\wt_\Gamma(x^{\lambda+1})=yz^{\nu_y},\quad\wt_\Gamma(z^{\nu+1})=x^\lambda,\quad\wt_\Gamma(yz^{\nu_y+1})=x^{\lambda-1}.
\]
\item If $\nu-\nu_y\geq t$, then $\nu=\nu_x+\nu_y+1$, $\lambda=2(t-\nu_y)-1$ and
\[
\wt_\Gamma(x^{\lambda+1})=z^\nu,\quad\wt_\Gamma(z^{\nu+1})=x^{\lambda-1},\quad\wt_\Gamma(yz^{\nu_y+1})=x^\lambda.
\]
\end{enumerate}
The symmetric statement holds if $z^{t-1},x,y^2\in\Gamma$.
\end{lem}

\begin{proof}
We can write $\Gamma=\Span(x^\lambda,xz^{\nu_x},yz^{\nu_y},z^\nu)$ with $\nu_x,\nu_y<\nu$.
Since $\wt(xz^t)=\wt(y)$ and $\wt(yz^t)=\wt(x)$, we have $\nu_x,\nu_y\leq t-1$. Since $\wt(z^{2t-1})=\wt(x^2)$, we have $\nu\leq2t-2$.
Write $\nu_y=t-m$ and $\nu=2t-n$ with $1\leq m\leq t$ and $2\leq n\leq t+1$. Note that $\nu-\nu_y=t+m-n$, so $\nu-\nu_y<t$ if and only if $m<n$.
Since $\wt(yz^{t-m})=\wt(x^{2m+1})$ and $\wt(z^{2t-n})=\wt(x^{2n})$, we have $x^{2m+1},x^{2n}\notin\Gamma$, so
\[
\lambda\leq
\begin{cases*}
2m & if $m<n$, \\
2n-1 & if $m\geq n$.
\end{cases*}
\]
Counting the monomials $z^q$ ($0\leq q\leq\nu$), $xz^q$ ($0\leq q\leq\nu_x$), $yz^q$ ($0\leq q\leq\nu_y$) and $x^i$ ($2\leq i\leq\lambda$), and using $\nu_x\leq t-1$, we obtain
\[
4t=(\nu+1)+(\nu_x+1)+(\nu_y+1)+(\lambda-1)\leq
\begin{cases*}
4t+m-n+1 & if $m<n$, \\
4t+n-m & if $m\geq n$,
\end{cases*}
\]
and in either case the right-hand side is at most $4t$.
Hence equality holds throughout: $\nu_x=t-1$, and
\[
\lambda=2m,\ n=m+1\ \text{ if }m<n,\qquad\lambda=2n-1,\ n=m\ \text{ if }m\geq n.
\]
In particular, since $\nu>\nu_x=t-1$, we have $n\leq t$, and
\[
\nu=2t-n=
\begin{cases*}
2t-m-1=\nu_x+\nu_y & if $m<n$, \\
2t-m=\nu_x+\nu_y+1 & if $m\geq n$.
\end{cases*}
\]

Assume that $m<n$, so that $\lambda=2m$ and $n=m+1$.
Then
\begin{align*}
\wt(x^{\lambda+1})&=\wt(x^{2m+1})=\wt(yz^{t-m})=\wt(yz^{\nu_y}), \\
\wt(z^{\nu+1})&=\wt(z^{2t-m})=\wt(x^{2m})=\wt(x^\lambda), \\
\wt(yz^{\nu_y+1})&=\wt(yz^{t-m+1})=\wt(x^{2m-1})=\wt(x^{\lambda-1}),
\end{align*}
and since $yz^{\nu_y},x^\lambda,x^{\lambda-1}\in\Gamma$, we obtain (i).

Assume that $m\geq n$, so that $\lambda=2n-1$ and $m=n$.
Then
\begin{align*}
\wt(x^{\lambda+1})&=\wt(x^{2n})=\wt(z^{2t-n})=\wt(z^\nu), \\
\wt(z^{\nu+1})&=\wt(z^{2t-n+1})=\wt(x^{2n-2})=\wt(x^{\lambda-1}), \\
\wt(yz^{\nu_y+1})&=\wt(yz^{t-n+1})=\wt(x^{2n-1})=\wt(x^\lambda),
\end{align*}
and since $z^\nu,x^{\lambda-1},x^\lambda\in\Gamma$, we obtain (ii).
\end{proof}

We now compute $k[S(\Gamma)]$ for the $G$-graphs classified in Lemmas \ref{z2l-1}, \ref{xzyz}, \ref{x2lyy2lx} and \ref{3 valleys}.

\begin{thm}\label{z2l-1 chart}
Let $\Gamma$ be one of the $G$-graphs appearing in Lemma \ref{z2l-1}, and assume that $\dim\sigma(\Gamma)=3$.
If $\Gamma=\Gamma_1$ or $\Gamma=\Gamma_2$, then $k[S(\Gamma)]\cong k[x,y,z]$. If $\Gamma=\Gamma_3$, then $k[S(\Gamma)]\cong k[x,y,z,w]/(xy-zw)$.
\end{thm}

\begin{proof}
By symmetry, we may assume that $\Gamma=\Gamma_1$ or $\Gamma=\Gamma_3$.

Assume $\Gamma=\Gamma_1=\Span(xz^{2t-1})$, so that $\lambda=1,\mu=0,\nu=2t-1$.
By Proposition \ref{monoid generator}, $S(\Gamma)$ is generated by $s_\Gamma(x^2)$, $s_\Gamma(y)$, $s_\Gamma(z^{2t})$.
Since $\wt_\Gamma(x^2)=z^{2t-1}$, $\wt_\Gamma(y)=xz^t$ and $\wt_\Gamma(z^{2t})=1$, we have
\[
s_\Gamma(x^2)=x^2z^{1-2t},\quad s_\Gamma(y)=x^{-1}yz^{-t},\quad s_\Gamma(z^{2t})=z^{2t},
\]
and Lemma \ref{basis} applies.

Now assume $\Gamma=\Gamma_3=\Span(xz^{t-1},z^{2t-1},yz^{t-1})$, so that $\lambda=\mu=1$, $\nu=2t-1$.
By Proposition \ref{monoid generator}, $S(\Gamma)$ is generated by
\[
s_\Gamma(x^2),\ s_\Gamma(y^2),\ s_\Gamma(z^{2t}),\ s_\Gamma(xy),\ s_\Gamma(xz^t),\ s_\Gamma(yz^t).
\]
Since
\begin{gather*}
\wt_\Gamma(x^2)=\wt_\Gamma(y^2)=z^{2t-1},\quad\wt_\Gamma(xz^t)=y,\\
\wt_\Gamma(yz^t)=x,\quad\wt_\Gamma(xy)=z^{t-1},\quad\wt_\Gamma(z^{2t})=1,
\end{gather*}
we have
\begin{align*}
s_\Gamma(x^2)&=x^2z^{1-2t}, & s_\Gamma(y^2)&=y^2z^{1-2t}, & s_\Gamma(z^{2t})&=z^{2t}, \\
s_\Gamma(xy)&=xyz^{1-t}, & s_\Gamma(xz^t)&=xy^{-1}z^t, & s_\Gamma(yz^t)&=x^{-1}yz^t.
\end{align*}
Since
\[
s_\Gamma(xy)=s_\Gamma(x^2)s_\Gamma(yz^t)=s_\Gamma(y^2)s_\Gamma(xz^t),\qquad
s_\Gamma(z^{2t})=s_\Gamma(xz^t)s_\Gamma(yz^t),
\]
$S(\Gamma)$ is generated by $s_\Gamma(x^2)$, $s_\Gamma(y^2)$, $s_\Gamma(xz^t)$, $s_\Gamma(yz^t)$, which satisfy $s_\Gamma(x^2)s_\Gamma(yz^t)=s_\Gamma(y^2)s_\Gamma(xz^t)$, and Lemma \ref{basis} applies.
\end{proof}

\begin{thm}\label{xzyz chart}
Let $\Gamma$ be one of the $G$-graphs appearing in Lemma \ref{xzyz}, and assume that $\dim\sigma(\Gamma)=3$.
Then $k[S(\Gamma)]\cong k[x,y,z]$.
\end{thm}

\begin{proof}
By symmetry, we may assume that $\Gamma=\Gamma_s^{xz}$.
By Proposition \ref{monoid generator}, $S(\Gamma)$ is generated by $s_\Gamma(x^{2s})$, $s_\Gamma(y)$, $s_\Gamma(z^{2t-s+1})$, $s_\Gamma(x^2z)$.
We have $\wt_\Gamma(x^{2s})=z^{2t-s}$, $\wt_\Gamma(z^{2t-s+1})=x^{2s-2}$ and $\wt_\Gamma(x^2z)=1$, while $\wt_\Gamma(y)=xz^t$ if $s\leq t$ and $\wt_\Gamma(y)=x^{2t+1}$ if $s=t+1$.
Hence, if $2\leq s\leq t$,
\begin{gather*}
s_\Gamma(x^{2s})=x^{2s}z^{s-2t},\quad s_\Gamma(y)=x^{-1}yz^{-t},\\
s_\Gamma(z^{2t-s+1})=x^{2-2s}z^{2t-s+1},\quad s_\Gamma(x^2z)=x^2z,
\end{gather*}
and if $s=t+1$,
\begin{gather*}
s_\Gamma(x^{2s})=x^{2t+2}z^{1-t},\quad s_\Gamma(y)=x^{-(2t+1)}y,\\
s_\Gamma(z^{2t-s+1})=x^{-2t}z^t,\quad s_\Gamma(x^2z)=x^2z.
\end{gather*}
In either case,
\[
s_\Gamma(x^2z)=s_\Gamma(x^{2s})\,s_\Gamma(z^{2t-s+1}),
\]
so $S(\Gamma)$ is generated by $s_\Gamma(x^{2s})$, $s_\Gamma(y)$, $s_\Gamma(z^{2t-s+1})$, and Lemma \ref{basis} applies.
\end{proof}

\begin{thm}\label{x2lyy2lx chart}
Let $\Gamma$ be one of the $G$-graphs appearing in Lemma \ref{x2lyy2lx}, and assume that $\dim\sigma(\Gamma)=3$.
Then $k[S(\Gamma)]\cong k[x,y,z]$.
\end{thm}

\begin{proof}
By symmetry, we may assume that $\Gamma=\Gamma_4$.
By Proposition \ref{monoid generator}, $S(\Gamma)$ is generated by $s_\Gamma(x^{2t+1})$, $s_\Gamma(z^t)$, $s_\Gamma(xy)$, $s_\Gamma(x^2z)$, $s_\Gamma(y^2)$, $s_\Gamma(yz)$.
Since
\begin{gather*}
\wt_\Gamma(x^{2t+1})=y,\quad\wt_\Gamma(z^t)=x^{2t},\quad\wt_\Gamma(xy)=z^{t-1},\\
\wt_\Gamma(x^2z)=1,\quad\wt_\Gamma(y^2)=x^2,\quad\wt_\Gamma(yz)=x^{2t-1},
\end{gather*}
we have
\begin{align*}
s_\Gamma(x^{2t+1})&=x^{2t+1}y^{-1}, & s_\Gamma(z^t)&=x^{-2t}z^t, & s_\Gamma(xy)&=xyz^{1-t}, \\
s_\Gamma(x^2z)&=x^2z, & s_\Gamma(y^2)&=x^{-2}y^2, & s_\Gamma(yz)&=x^{1-2t}yz.
\end{align*}
Since
\begin{gather*}
s_\Gamma(yz)=s_\Gamma(z^t)s_\Gamma(xy),\quad
s_\Gamma(x^2z)=s_\Gamma(x^{2t+1})s_\Gamma(yz),\\
s_\Gamma(y^2)=s_\Gamma(x^{2t+1})^{t-2}s_\Gamma(z^t)^{t-1}s_\Gamma(xy)^t,
\end{gather*}
$S(\Gamma)$ is generated by $s_\Gamma(x^{2t+1})$, $s_\Gamma(z^t)$, $s_\Gamma(xy)$, and Lemma \ref{basis} applies.
\end{proof}

\begin{thm}\label{3 valleys chart}
Let $\Gamma$ be a $G$-graph as in Lemma \ref{3 valleys}, and assume that $\dim\sigma(\Gamma)=3$.
Then $k[S(\Gamma)]\cong k[x,y,z]$.
\end{thm}

\begin{proof}
By symmetry, we may assume that $x^2,y\in\Gamma$.
By Proposition \ref{monoid generator}, $S(\Gamma)$ is generated by $s_\Gamma(x^{\lambda+1})$, $s_\Gamma(z^{\nu+1})$, $s_\Gamma(xy)$, $s_\Gamma(x^2z)$, $s_\Gamma(y^2)$, $s_\Gamma(xz^t)$, $s_\Gamma(yz^{\nu_y+1})$.
Write $m\coloneqq t-\nu_y\geq1$.
We have $\wt_\Gamma(xz^t)=y$, $\wt_\Gamma(x^2z)=1$, $\wt_\Gamma(xy)=z^{t-1}$ and $\wt_\Gamma(y^2)=x^2$.

Suppose first that $\nu-\nu_y<t$.
By Lemma \ref{3 valleys}, $\lambda=2m$, $\nu=\nu_x+\nu_y$, $\wt_\Gamma(x^{\lambda+1})=yz^{\nu_y}$, $\wt_\Gamma(z^{\nu+1})=x^\lambda$, $\wt_\Gamma(yz^{\nu_y+1})=x^{\lambda-1}$.
Hence
\begin{align*}
s_\Gamma(x^{\lambda+1})&=x^{\lambda+1}y^{-1}z^{-\nu_y}, & s_\Gamma(z^{\nu+1})&=x^{-\lambda}z^{\nu+1}, & s_\Gamma(xy)&=xyz^{1-t}, \\
s_\Gamma(x^2z)&=x^2z, & s_\Gamma(y^2)&=x^{-2}y^2, & s_\Gamma(xz^t)&=xy^{-1}z^t, \\
s_\Gamma(yz^{\nu_y+1})&=x^{1-\lambda}yz^{\nu_y+1}.
\end{align*}
One checks directly that
\begin{gather*}
s_\Gamma(xz^t)=s_\Gamma(x^{\lambda+1})s_\Gamma(z^{\nu+1}),\quad
s_\Gamma(x^2z)=s_\Gamma(xz^t)s_\Gamma(xy),\\
s_\Gamma(yz^{\nu_y+1})=s_\Gamma(xy)s_\Gamma(z^{\nu+1}).
\end{gather*}
If $m\geq2$, then also
\[
s_\Gamma(y^2)=s_\Gamma(x^{\lambda+1})^{m-2}s_\Gamma(z^{\nu+1})^{m-1}s_\Gamma(xy)^m,
\]
so $S(\Gamma)$ is generated by $s_\Gamma(x^{\lambda+1})$, $s_\Gamma(z^{\nu+1})$, $s_\Gamma(xy)$.
If $m=1$ (so $\lambda=2$, $\nu_y=t-1$), then instead
\[
s_\Gamma(xy)=s_\Gamma(x^{\lambda+1})s_\Gamma(y^2),
\]
so $S(\Gamma)$ is generated by $s_\Gamma(x^{\lambda+1})$, $s_\Gamma(z^{\nu+1})$, $s_\Gamma(y^2)$.
In either case Lemma \ref{basis} applies.

Now suppose $\nu-\nu_y\geq t$.
By Lemma \ref{3 valleys}, $\lambda=2m-1$, $\nu=\nu_x+\nu_y+1$, $\wt_\Gamma(x^{\lambda+1})=z^\nu$, $\wt_\Gamma(z^{\nu+1})=x^{\lambda-1}$, $\wt_\Gamma(yz^{\nu_y+1})=x^\lambda$. Note that $m\geq2$, since $\lambda\geq2$.
Hence
\[
s_\Gamma(x^{\lambda+1})=x^{\lambda+1}z^{-\nu},\quad
s_\Gamma(z^{\nu+1})=x^{1-\lambda}z^{\nu+1},\quad
s_\Gamma(yz^{\nu_y+1})=x^{-\lambda}yz^{\nu_y+1}.
\]
One checks directly that
\begin{gather*}
s_\Gamma(x^2z)=s_\Gamma(x^{\lambda+1})s_\Gamma(z^{\nu+1}),\quad
s_\Gamma(xy)=s_\Gamma(x^{\lambda+1})s_\Gamma(yz^{\nu_y+1}),\\
s_\Gamma(z^{\nu+1})=s_\Gamma(xz^t)s_\Gamma(yz^{\nu_y+1}),
\end{gather*}
and, since $m\geq2$,
\[
s_\Gamma(y^2)=s_\Gamma(x^{\lambda+1})^{m-1}s_\Gamma(xz^t)^{m-2}s_\Gamma(yz^{\nu_y+1})^m.
\]
Thus $S(\Gamma)$ is generated by $s_\Gamma(x^{\lambda+1})$, $s_\Gamma(xz^t)$, $s_\Gamma(yz^{\nu_y+1})$, and Lemma \ref{basis} applies.
\end{proof}

\subsection{\texorpdfstring{The case $z^{t-1}\notin\Gamma$}{The case z\string^(t-1) ∉ Γ}}\label{subsec:z notin}

We treat separately the cases where $x,y\in\Gamma$ and where $x\notin\Gamma$ or $y\notin\Gamma$.

\begin{lem}\label{level structure}
Let $\Gamma$ be a $G$-graph with $z^{t-1}\notin\Gamma$ such that $x,y\in\Gamma$ (so $\nu\leq t-2$).
Then exactly one of $x^2,y^2$ lies in $\Gamma$.
Moreover, for every $1\leq q\leq\nu$, at least one of $xz^q,yz^q$ lies in $\Gamma$, and if both do, then so does $xyz^q$.
\end{lem}

\begin{proof}
Since $\wt(x^2)=\wt(y^2)$, the monomials $x^2,y^2$ do not both lie in $\Gamma$.
If neither does, then every monomial in $\Gamma$ is of the form $x^{\delta_1}y^{\delta_2}z^q$ with $\delta_1,\delta_2\in\{0,1\}$ and $0\leq q\leq\nu$, so $|\Gamma|\leq4(\nu+1)\leq4(t-1)<4t$, a contradiction.
Hence exactly one of $x^2,y^2$ lies in $\Gamma$. By symmetry, we may assume that $x^2\in\Gamma$, so that $\mu=1$.

Suppose $xz^q,yz^q\in\Gamma$ for some $0\leq q\leq\nu$, and set $u\coloneqq xyz^q$.
If $u\notin\Gamma$, then, since $u/x=yz^q\in\Gamma$ and $u/y=xz^q\in\Gamma$, Lemma \ref{not divisible} shows that $\wt_\Gamma(u)$ is divisible by neither $x$ nor $y$, so $\wt_\Gamma(u)=z^j$ for some $j\geq0$.
Since $\wt(xy)=\wt(z^{t-1})$, we have $\wt(z^j)=\wt(u)=\wt(z^{t-1+q})$.
Since $\wt(z)=\wt(x^{-2})$ and $\wt(x)$ has order $4t$, $\wt(z)$ has order $2t$, so $j\equiv t-1+q\pmod{2t}$.
As $0\leq q\leq t-2$, we have $t-1\leq t-1+q\leq2t-3<2t$, so $j\geq t-1+q\geq t-1$.
Hence $z^{t-1}\in\Gamma$, a contradiction.
Hence $xyz^q\in\Gamma$.

Set $c\coloneqq\max\{0\leq q\leq\nu\mid xz^q\in\Gamma\}$, $b\coloneqq\max\{0\leq q\leq\nu\mid yz^q\in\Gamma\}$ and $\beta\coloneqq\max\{i\geq0\mid x^iy\in\Gamma\}$.
We count $|\Gamma|$.
Since $\wt(x^2z)=\wt(y^2z)=\wt(1)$ and $\mu=1$, the previous paragraph shows that the monomials in $\Gamma$ are the following: $x^i$ ($0\leq i\leq\lambda$), $x^iy$ ($0\leq i\leq\beta$), and, for $1\leq q\leq\nu$, the monomials $z^q$, $xz^q$ ($q\leq c$), $yz^q$ ($q\leq b$) and $xyz^q$ ($q\leq\min\{c,b\}$).
Since $\wt(yz^b)=\wt(x^{2(t-b)+1})$ and $yz^b\in\Gamma$, we have $\lambda\leq2t-2b$.
Since $\wt(z^\nu)=\wt(x^{2(t-1-\nu)+1}y)$ and $z^\nu\in\Gamma$, we have $\beta\leq2(t-1-\nu)$.
Hence
\begin{equation}\label{level count}
\begin{aligned}
4t=|\Gamma|&=(\lambda+1)+(\beta+1)+\nu+c+b+\min\{c,b\} \\
&\leq(2t-2b+1)+(2t-2\nu-1)+\nu+c+b+\min\{c,b\} \\
&=4t+c-\nu-(b-\min\{c,b\}).
\end{aligned}
\end{equation}

Finally, suppose $xz^q,yz^q\notin\Gamma$ for some $1\leq q\leq\nu$.
Then $c,b<q\leq\nu$, and \eqref{level count} gives $4t\leq4t+c-\nu<4t$, a contradiction.
\end{proof}

We can now determine the shape of $\Gamma$.

\begin{cor}\label{level structure shape}
Let $\Gamma$ be as in Lemma \ref{level structure}, and assume by symmetry that $x^2\in\Gamma$.
Set $b\coloneqq\max\{0\leq q\leq\nu\mid yz^q\in\Gamma\}$ and $\beta\coloneqq\max\{i\geq0\mid x^iy\in\Gamma\}$.
Then $xz^\nu\in\Gamma$,
\[
\lambda=2t-2b,\qquad \mu=1,\qquad \beta=2(t-1-\nu),
\]
and
\[
\Gamma=\Span(x^\lambda,x^\beta y,xyz^b,xz^\nu).
\]
\end{cor}

\begin{proof}
Since $y\in\Gamma$ and $y^2\notin\Gamma$, we have $\mu=1$.
Let $c$ be as in the proof of Lemma \ref{level structure}.
By Lemma \ref{level structure}, $c=\nu$ or $b=\nu$.
If $c<\nu$, then $b=\nu$, and \eqref{level count} gives $4t\leq4t+2c-2\nu<4t$, a contradiction.
Hence $c=\nu$, that is, $xz^\nu\in\Gamma$. Then $\min\{c,b\}=b$, so equality holds in \eqref{level count}, which gives $\lambda=2t-2b$ and $\beta=2(t-1-\nu)$.
By the description of the monomials in $\Gamma$ in the proof of Lemma \ref{level structure}, we obtain $\Gamma=\Span(x^\lambda,x^\beta y,xyz^b,xz^\nu)$.
\end{proof}

\begin{thm}\label{level structure chart}
Let $\Gamma$ be as in Corollary \ref{level structure shape}, and assume that $\dim\sigma(\Gamma)=3$.
If $b=\nu$, then $k[S(\Gamma)]\cong k[x,y,z]$.
If $b<\nu$, then $k[S(\Gamma)]\cong k[x,y,z,w]/(xy-zw)$.
\end{thm}

\begin{proof}
By Proposition \ref{monoid generator}, the monoid $S(\Gamma)$ is generated by $s_\Gamma(x^{\lambda+1})$, $s_\Gamma(y^2)$, $s_\Gamma(z^{\nu+1})$, $s_\Gamma(x^2z)$, $s_\Gamma(x^{\beta+1}y)$, together with $s_\Gamma(yz^{b+1})$ if $b<\nu$.
Recall that $\lambda=2t-2b$ and $\beta=2(t-1-\nu)\geq2$.
By Remark \ref{wt relations},
\begin{align*}
\wt_\Gamma(x^{\lambda+1})&=yz^b, & \wt_\Gamma(y^2)&=x^2, & \wt_\Gamma(z^{\nu+1})&=x^{\beta-1}y, \\
\wt_\Gamma(x^{\beta+1}y)&=z^\nu, & \wt_\Gamma(x^2z)&=1, & \wt_\Gamma(yz^{b+1})&=x^{\lambda-1}\ (\text{if }b<\nu),
\end{align*}
and hence
\begin{align*}
s_\Gamma(x^{\lambda+1})&=x^{\lambda+1}y^{-1}z^{-b}, & s_\Gamma(y^2)&=x^{-2}y^2, & s_\Gamma(z^{\nu+1})&=x^{1-\beta}y^{-1}z^{\nu+1}, \\
s_\Gamma(x^2z)&=x^2z, & s_\Gamma(x^{\beta+1}y)&=x^{\beta+1}yz^{-\nu}, & s_\Gamma(yz^{b+1})&=x^{1-\lambda}yz^{b+1}.
\end{align*}
In either case,
\[
s_\Gamma(x^2z)=s_\Gamma(z^{\nu+1})\,s_\Gamma(x^{\beta+1}y).
\]

Suppose $b=\nu$.
Then $\lambda-1=\beta+1$, so $s_\Gamma(x^{\beta+1}y)=s_\Gamma(x^{\lambda+1})\,s_\Gamma(y^2)$, and hence $S(\Gamma)$ is generated by $s_\Gamma(x^{\lambda+1})$, $s_\Gamma(y^2)$, $s_\Gamma(z^{\nu+1})$, and Lemma \ref{basis} applies.

Suppose $b<\nu$.
Then $s_\Gamma(x^2z)=s_\Gamma(x^{\lambda+1})\,s_\Gamma(yz^{b+1})$, and (as $\nu-b\geq1$)
\[
s_\Gamma(y^2)=s_\Gamma(x^{\beta+1}y)\,s_\Gamma(yz^{b+1})^{\nu-b}\,s_\Gamma(x^{\lambda+1})^{\nu-b-1},
\]
so $S(\Gamma)$ is generated by $s_\Gamma(x^{\lambda+1}),s_\Gamma(x^{\beta+1}y),s_\Gamma(z^{\nu+1}),s_\Gamma(yz^{b+1})$, which satisfy $s_\Gamma(z^{\nu+1})s_\Gamma(x^{\beta+1}y)=s_\Gamma(x^{\lambda+1})s_\Gamma(yz^{b+1})$, and Lemma \ref{basis} applies.
\end{proof}

We now turn to the case where $x\notin\Gamma$ or $y\notin\Gamma$.

\begin{lem}\label{y absent shape}
Let $\Gamma$ be a $G$-graph with $z^{t-1}\notin\Gamma$ and $y\notin\Gamma$.
Then $x\in\Gamma$ and
\[
\lambda=4t-1-2\nu,\qquad \Gamma=\Span(x^\lambda,xz^\nu).
\]
The symmetric statement holds if $z^{t-1}\notin\Gamma$ and $x\notin\Gamma$.
\end{lem}

\begin{proof}
If $x\notin\Gamma$, then $\Gamma\subset\{z^q\mid0\leq q\leq t-2\}$, so $|\Gamma|\leq t-1<4t$, a contradiction. Hence $x\in\Gamma$.
Suppose $xz^\nu\notin\Gamma$.
Since $z^\nu\in\Gamma$, Lemma \ref{not divisible} shows that $\wt_\Gamma(xz^\nu)$ is not divisible by $x$. Since $y\notin\Gamma$, no element of $\Gamma$ is divisible by $y$, so $\wt_\Gamma(xz^\nu)=z^j$ for some $j\geq0$.
Then $\wt(xz^\nu)=\wt(z^j)$, that is, $1-2\nu\equiv-2j\pmod{4t}$. This is impossible, since the left-hand side is odd while the right-hand side and $4t$ are even.
Hence $xz^\nu\in\Gamma$.

Since $\Gamma$ is closed under taking divisors, $xz^q\in\Gamma$ for $0\leq q\leq\nu$, and $x^2z^q\notin\Gamma$ for $1\leq q\leq\nu$, so $\Gamma\cap\{x^iz^q\mid q\geq1\}=\{z^q,xz^q\mid1\leq q\leq\nu\}$.
Hence $\Gamma=\Span(x^\lambda,xz^\nu)$, and $4t=|\Gamma|=(\lambda+1)+2\nu$ gives $\lambda=4t-1-2\nu$.
\end{proof}

\begin{thm}\label{y absent chart}
Let $\Gamma$ be as in Lemma \ref{y absent shape}, and assume that $\dim\sigma(\Gamma)=3$.
Then $k[S(\Gamma)]\cong k[x,y,z]$.
\end{thm}

\begin{proof}
By Proposition \ref{monoid generator}, $S(\Gamma)$ is generated by $s_\Gamma(x^{\lambda+1})$, $s_\Gamma(y)$, $s_\Gamma(z^{\nu+1})$, and, if $\nu\geq1$, $s_\Gamma(x^2z)$.
Since $\lambda=4t-1-2\nu\geq2t+3$, we have $\wt_\Gamma(x^{\lambda+1})=z^\nu$, $\wt_\Gamma(y)=x^{2t+1}$ and $\wt_\Gamma(z^{\nu+1})=x^{\lambda-1}$, so
\begin{align*}
s_\Gamma(x^{\lambda+1})&=x^{\lambda+1}z^{-\nu}, & s_\Gamma(y)&=x^{-(2t+1)}y, \\
s_\Gamma(z^{\nu+1})&=x^{1-\lambda}z^{\nu+1}.
\end{align*}
If $\nu\geq1$, then $s_\Gamma(x^2z)=x^2z=s_\Gamma(x^{\lambda+1})s_\Gamma(z^{\nu+1})$, so in either case $S(\Gamma)$ is generated by $s_\Gamma(x^{\lambda+1}),s_\Gamma(y),s_\Gamma(z^{\nu+1})$, and Lemma \ref{basis} applies.
\end{proof}

Combining the results of Sections \ref{subsec:z in} and \ref{subsec:z notin}, we obtain Theorem \ref{thm:intro conifold} for groups of type (4).

\begin{cor}\label{Hilb conifold section4}
Let $G$ be of type $\frac{1}{4t}(1,2t+1,4t-2)$ with $t\geq2$, and let $\Gamma$ be a $G$-graph with $\dim\sigma(\Gamma)=3$.
Then $V(\Gamma)=U_{\sigma(\Gamma)}$.
Moreover, $\Hilb^G(\A_k^3)$ has at worst conifold singularities.
In particular, $\Hilb^G(\A_k^3)$ is normal.
\end{cor}

\begin{proof}
We check that every $G$-graph falls under one of Theorems \ref{z2l-1 chart}, \ref{xzyz chart}, \ref{x2lyy2lx chart}, \ref{3 valleys chart}, \ref{level structure chart} and \ref{y absent chart}.
Suppose first that $z^{t-1}\in\Gamma$.
If $x^2,y^2\notin\Gamma$, we are in Lemma \ref{z2l-1}.
Otherwise, by symmetry, we may assume that $x^2\in\Gamma$. If $y\notin\Gamma$, we are in Lemma \ref{xzyz}, and if $y\in\Gamma$, then by Lemmas \ref{xxyyyx} and \ref{x2lyy2lx}, $\Gamma$ has two or three valleys, so we are in Lemma \ref{x2lyy2lx} or Lemma \ref{3 valleys}.
Suppose next that $z^{t-1}\notin\Gamma$.
If $x,y\in\Gamma$, we are in Corollary \ref{level structure shape}. Otherwise, we are in Lemma \ref{y absent shape}, up to symmetry.
Hence $k[S(\Gamma)]$ is isomorphic to $k[x,y,z]$ or $k[x,y,z,w]/(xy-zw)$.
Both are normal, so $V(\Gamma)=U_{\sigma(\Gamma)}$.
The remaining assertions follow from Theorem \ref{Nakamura fan}.
\end{proof}

\section{The general canonical case}\label{sec:FB}

In this section, we prove Theorem \ref{thm:intro canonical}.
Assertions (i) and (ii) follow by combining Sections \ref{sec:type3} and \ref{sec:type4} with the known cases (Section \ref{subsec:Hilb canonical}).
For (iii), we first recall universal flattenings and F-blowups.
We then show that for $l\geq2$, the F-blowup $\FB_{(l)}$ of the conifold singularity is its blowup at the singular point (Proposition \ref{conifold resolution}), and apply this to the charts of $\Hilb^G(\A_k^3)$ (Corollary \ref{FB FB resolution}).

\subsection{\texorpdfstring{$G$-Hilbert schemes}{G-Hilbert schemes}}\label{subsec:Hilb canonical}

Combining Sections \ref{sec:type3} and \ref{sec:type4} with the known cases, we obtain the following.

\begin{thm}\label{Hilb canonical}
Let $G\subset\GL(3,k)$ be a finite abelian subgroup containing no pseudo-reflections, of order prime to the characteristic of $k$, such that $X=\A_k^3/G$ has at worst a canonical singularity.
Then:
\begin{enumerate}
\item $\Hilb^G(\A_k^3)$ has at worst conifold singularities; in particular, it is normal;
\item every essential divisor over $X$ appears on $\Hilb^G(\A_k^3)$ as a prime divisor.
\end{enumerate}
\end{thm}

\begin{proof}
By Proposition \ref{canonical classification} and the subsequent remark, either $G\subset\SL(3,k)$ or $G$ is of one of the types (2)--(5).
If $G\subset\SL(3,k)$, then $\Hilb^G(\A_k^3)$ is smooth by Nakamura \cite{Nak01} (see also \cite{BKR01}). It has at worst conifold singularities for type (2) by K\k{e}dzierski \cite{Ked11}, and for types (3) and (4) by Corollaries \ref{Hilb conifold} and \ref{Hilb conifold section4}. For type (5), it is smooth by Remark \ref{type5 smooth}.
This proves (i), as conifold singularities are normal.
Since the Hilbert--Chow morphism $\Hilb^G(\A_k^3)\to X$ is proper, birational and toric, (ii) follows from (i) and Lemma \ref{conifold ess}.
\end{proof}

\subsection{Universal flattening}

We recall the universal flattening of a coherent sheaf.
For a direct sum of fractional ideals, we reduce the flatness condition to the invertibility of a single ideal (Lemmas \ref{flattening tor} and \ref{product invertible}).

Throughout this subsection, $X$ denotes a Noetherian integral scheme with function field $K(X)$, and $\mathcal{K}=\mathcal{K}_X$ denotes the constant sheaf on $X$ associated to $K(X)$.
For a birational morphism $\pi\colon Y\to X$, we identify $\pi^\ast\mathcal{K}_X$ with $\mathcal{K}_Y$.

\begin{dfn}
For a coherent sheaf $\mathcal F$ of $\mathcal O_X$-modules, we define the \emph{torsion subsheaf} $\tor(\mathcal F)\subset\mathcal F$ to be the kernel of the natural map $\mathcal F\to\mathcal F\otimes_{\mathcal O_X}\mathcal K_X$.
\end{dfn}

\begin{dfn}[\cite{Vil06}]\label{univ flattening}
Let $\mathcal F$ be a coherent sheaf of $\mathcal O_X$-modules.
We say that a proper birational morphism $\pi\colon Y\to X$ is the \emph{universal flattening} of $\mathcal F$ if
\begin{enumerate}
\item $\pi^\ast\mathcal F/\tor(\pi^\ast\mathcal F)$ is locally free, and
\item for any proper birational morphism $\pi'\colon Y'\to X$ such that the sheaf $(\pi')^\ast\mathcal F/\tor((\pi')^\ast\mathcal F)$ is locally free, there exists a unique morphism $\rho\colon Y'\to Y$ with $\pi'=\pi\circ\rho$.
\end{enumerate}
\end{dfn}

\begin{rem}
The universal flattening of $\mathcal F$ exists and is unique up to a unique isomorphism over $X$. It may be realized as a blowup along a suitable fractional ideal attached to $\mathcal F$ \cite[Theorem 3.3]{Vil06}.
\end{rem}

\begin{lem}\label{flattening tor}
Let $\mathcal I_1,\ldots,\mathcal I_n\subset\mathcal K$ be coherent fractional ideal sheaves, that is, coherent $\mathcal O_X$-submodules of $\mathcal K$ with $\mathcal I_j\otimes_{\mathcal O_X}\mathcal K=\mathcal K$ for all $j$, and let $\pi\colon Y\to X$ be a birational morphism.
Then the following are equivalent:
\begin{enumerate}
\item $\pi^\ast\big(\bigoplus_j\mathcal I_j\big)\big/\tor\big(\pi^\ast\big(\bigoplus_j\mathcal I_j\big)\big)$ is locally free.
\item $\mathcal I_j\cdot\mathcal{O}_Y$ is invertible for all $j$.
\end{enumerate}
\end{lem}

\begin{proof}
Since $\pi^\ast\mathcal K=\mathcal K_Y$, the natural map $\pi^\ast\mathcal I_j\to\pi^\ast\mathcal I_j\otimes_{\mathcal O_Y}\mathcal K_Y=\mathcal K_Y$ has kernel $\tor(\pi^\ast\mathcal I_j)$ and image $\mathcal I_j\cdot\mathcal O_Y$. Since everything commutes with finite direct sums,
\[
\pi^\ast\Big(\bigoplus_j\mathcal I_j\Big)\Big/\tor\Big(\pi^\ast\Big(\bigoplus_j\mathcal I_j\Big)\Big)\cong\bigoplus_j\big(\mathcal I_j\cdot\mathcal O_Y\big).
\]
Since a direct summand of a free module over a local ring is free, this is locally free if and only if each summand is locally free.
Since each summand has generic rank $1$, this means that each summand is invertible.
\end{proof}

\begin{lem}\label{product invertible}
Let $\mathcal I_1,\ldots,\mathcal I_n\subset\mathcal K$ be coherent fractional ideal sheaves and let $\pi\colon Y\to X$ be a birational morphism.
Then $(\mathcal I_1\cdots\mathcal I_n)\cdot\mathcal O_Y$ is invertible if and only if $\mathcal I_j\cdot\mathcal O_Y$ is invertible for all $j$.
\end{lem}

\begin{proof}
The ``if'' part is clear.
Conversely, fix $y\in Y$ and let $\mathfrak a_j\coloneqq\mathcal I_j\cdot\mathcal O_{Y,y}$.
If $\mathfrak a_1\cdots\mathfrak a_n=t\mathcal O_{Y,y}$ with $t\in K(X)^\times$, then $\mathfrak a_j\cdot t^{-1}\prod_{i\neq j}\mathfrak a_i=\mathcal O_{Y,y}$, so each $\mathfrak a_j$ is an invertible fractional ideal of the local domain $\mathcal O_{Y,y}$, hence principal.
\end{proof}

\begin{cor}\label{flattening multiplicity}
Let $\mathcal F=\bigoplus_\gamma\mathcal I_\gamma$ be a finite direct sum of coherent fractional ideal sheaves, and let $\mathcal J_1,\ldots,\mathcal J_s$ be the distinct sheaves occurring among the $\mathcal I_\gamma$, disregarding multiplicity.
Then the universal flattening of $\mathcal F$ coincides with the universal flattening of $\bigoplus_{i=1}^s\mathcal J_i$.
\end{cor}

\begin{proof}
By Lemma \ref{flattening tor}, for both sheaves the flatness condition on $\pi\colon Y\to X$ says that every $\mathcal J_i\cdot\mathcal O_Y$ is invertible.
\end{proof}

\subsection{F-blowups}

Let $X$ be a normal toric variety with lattice of characters $M$, and let $l$ be a positive integer.
For an affine toric open subset $U=\Spec k[A]\subset X$, where $A=\sigma^\vee\cap M$, the inclusion $A\hookrightarrow(1/l)A\coloneqq\sigma^\vee\cap(1/l)M$ induces a finite toric morphism $\Spec k[(1/l)A]\to U$. These glue to the \emph{Frobenius-like morphism} $F^{(l)}\colon X_{(l)}\to X$ \cite[Section 2.3]{CMDY25}.
If $k$ has characteristic $p>0$ and $l=p^e$, then $F^{(l)}$ is the $e$-th iterated Frobenius morphism.
For $X=\Spec R$ affine with $R=k[A]$, we write $R^{1/l}\coloneqq k[(1/l)A]=F^{(l)}_\ast\mathcal O_{X_{(l)}}$, regarded as an $R$-module.

\begin{dfn}[\cite{Yas12}]
We define the \emph{$(l)$-th F-blowup} of $X$, denoted by $\FB_{(l)}(X)$, to be the universal flattening of $F^{(l)}_\ast\mathcal{O}_{X_{(l)}}$ (Definition \ref{univ flattening}).
\end{dfn}

\begin{dfn}
For a normal toric variety $X$, the F-blowups $\FB_{(l)}(X)$ are isomorphic to each other for $l\gg0$ \cite[Theorem 3.13]{Yas12}. We define $\FB_{(\infty)}(X)\coloneqq\FB_{(l)}(X)$ for $l\gg0$.
\end{dfn}

\begin{rem}[{See \cite[Remark 2.5]{CMDY25}}]\label{Hilb FB iso}
Suppose that $G$ is a finite abelian group of order prime to the characteristic of $k$ acting linearly on $\A_{k}^d$, and let $X\coloneqq\A_{k}^d/G$.
We have an isomorphism $\Hilb^G(\A_{k}^d)\cong\FB_{(\infty)}(X)$ compatible with the birational morphism to $X$ \cite{Yas12,TY09}.
For the analogous statement for finite linearly reductive group schemes, such as $\boldsymbol{\mu}_p$ in characteristic $p$, see \cite[Theorem 4.3]{LY24}.
\end{rem}

\subsection{F-blowups of the conifold singularity}

By Theorem \ref{Hilb canonical}(i), the conifold singularity is the only singularity that appears on $\Hilb^G(\A_k^3)$.
We now compute its F-blowups.

\begin{prop}[{cf.\ \cite[Example 4.14 and Proposition 4.15]{FMS19}}]\label{conifold decomposition}
Let $R=k[x,y,z,w]/(xy-zw)$ be the coordinate ring of the conifold singularity.
Then for $l\geq1$, the $R$-module $R^{1/l}$ is a direct sum of modules isomorphic to $R$, $(x,z)$ or $(x,w)$.
Moreover, for every $l\geq2$, both $(x,z)$ and $(x,w)$ occur as direct summands.
\end{prop}

\begin{proof}
Let $M=\Z^3$, and let $C$ be the cone generated by $(0,0,1)$, $(1,0,1)$, $(1,1,1)$, $(0,1,1)$, which are the exponents of $x$, $z$, $y$, $w$, respectively.
Then $R=k[C\cap M]$ and $C=\{(m_1,m_2,m_3)\mid m_1,m_2\geq0,\ m_3\geq m_1,\ m_3\geq m_2\}$.
By \cite[Proposition 4.15]{FMS19}, $R^{1/l}\cong\bigoplus_{u\in M/lM}k[(C+u/l)\cap M]$, and by \cite[Example 4.14]{FMS19}, each summand is isomorphic to $R$, $(x,z)$ or $(x,w)$.
Let $l\geq2$.
For $u=(1,0,0)$, the lattice points of $C+u/l$ are those $m\in C$ with $m_1\geq1$, so the corresponding summand is the ideal $(y,z)=(y/w)\cdot(x,w)\cong(x,w)$. Similarly, $u=(0,1,0)$ gives $(y,w)=(y/z)\cdot(x,z)\cong(x,z)$.
\end{proof}

See also \cite[Section 6]{TY09} for a related statement for the complete local ring in odd characteristic.

\begin{prop}\label{conifold resolution}
Let $R=k[x,y,z,w]/(xy-zw)$ and let $X=\Spec R$ be the conifold singularity.
Then $\FB_{(l)}(X)$ is a resolution of $X$ for every $l\geq2$. In particular, $\FB_{(\infty)}(X)$ is a resolution of $X$.
\end{prop}

\begin{proof}
Fix $l\geq2$.
By Proposition \ref{conifold decomposition} and Corollary \ref{flattening multiplicity}, $\FB_{(l)}(X)$, the universal flattening of $R^{1/l}$, coincides with the universal flattening of $R\oplus(x,z)\oplus(x,w)$.

Let $\pi\colon Y\to X$ be a birational morphism.
Since $R\cdot\mathcal O_Y=\mathcal O_Y$ is always invertible, Lemma \ref{flattening tor} shows that $\pi^\ast(R\oplus(x,z)\oplus(x,w))/\tor(\cdots)$ is locally free if and only if $(x,z)\cdot\mathcal O_Y$ and $(x,w)\cdot\mathcal O_Y$ are both invertible.
By Lemma \ref{product invertible}, this holds if and only if $\big((x,z)(x,w)\big)\cdot\mathcal O_Y$ is invertible.
In $R$ we have $zw=xy$, so
\[
(x,z)(x,w)=(x^2,xw,xz,zw)=(x^2,xw,xz,xy)=x\cdot(x,y,z,w)=x\mathfrak m,
\]
where $\mathfrak m=(x,y,z,w)$ is the maximal ideal at the vertex of $X$.
Since $(x)$ is a principal ideal generated by a nonzerodivisor, it is invertible on any $Y$. Hence, by Lemma \ref{product invertible} applied to $\mathcal I_1=(x)$ and $\mathcal I_2=\mathfrak m$, $(x\mathfrak m)\cdot\mathcal O_Y$ is invertible if and only if $\mathfrak m\cdot\mathcal O_Y$ is invertible.

Hence $\FB_{(l)}(X)$ is the universal flattening of $\mathfrak m$, that is, the blowup of $X$ at $\mathfrak m$.
The origin $x=y=z=w=0$ is an ordinary double point on the hypersurface $X\subset\A_k^4$, and its projectivized tangent cone is the smooth quadric $\{xy=zw\}\cong\mathbb P^1\times\mathbb P^1\subset\mathbb P^3$. Therefore, by \cite[(1.9), Example (1), pp.~351--352]{Rei87}, the blowup of $X$ at this point (equivalently, at $\mathfrak m$) is smooth, with exceptional divisor isomorphic to $\mathbb P^1\times\mathbb P^1$.
\end{proof}

\begin{rem}
The exceptional divisor of $\FB_{(l)}(X)\to X$ ($l\geq2$) corresponds to the ray through the center of the square cone of $X$. This divisor is essential over $X$, although $X$ has no BGS essential exceptional divisor \cite[Example A.7]{CMDY25}.
\end{rem}

\subsection{\texorpdfstring{F-blowups of $G$-Hilbert schemes}{F-blowups of G-Hilbert schemes}}

We can now prove Theorem \ref{thm:intro canonical}(iii).

\begin{cor}\label{FB FB resolution}
Let $G$ and $X$ be as in Theorem \ref{Hilb canonical}.
Then $\FB_{(l)}(\Hilb^G(\A_k^3))$ is a resolution of $X$ for every $l\geq2$. In particular, so is $\FB_{(\infty)}(\Hilb^G(\A_k^3))$, which is isomorphic to $\FB_{(\infty)}(\FB_{(\infty)}(X))$.
\end{cor}

\begin{proof}
The isomorphism follows from Remark \ref{Hilb FB iso}.
Set $Y=\Hilb^G(\A_k^3)$ and fix $l\geq2$.
By Theorem \ref{Hilb canonical}, $Y$ is a normal toric variety covered by affine toric open subsets $U_\tau$, each isomorphic to $\A_k^3$ or to a conifold singularity.
Since $\FB_{(l)}(-)$ is compatible with open immersions \cite[Proposition 3.2(3)]{Yas12}, it suffices to check that each $\FB_{(l)}(U_\tau)$ is smooth.
If $U_\tau\cong\A_k^3$, then $F^{(l)}_\ast\mathcal O_{(U_\tau)_{(l)}}$ is free, so $\FB_{(l)}(U_\tau)=U_\tau$. Otherwise, use Proposition \ref{conifold resolution}.
\end{proof}

\begin{rem}\label{FB conifold}
Theorem \ref{Hilb canonical} and Corollary \ref{FB FB resolution} hold for $\FB_{(\infty)}(X)$ for every three-dimensional $\Q$-factorial canonical toric singularity $X$, in any characteristic.
Indeed, $\FB_{(\infty)}(X)$ is determined by the defining cone of $X$, independently of $k$ \cite[Proposition 3.12 and Theorem 3.13]{Yas12}.
Hence we may assume that the characteristic of $k$ does not divide $|N/N_0|$, where $N_0\subset N$ is the sublattice generated by the primitive generators of the rays. Then $X\cong\A_k^3/G$ and $\FB_{(\infty)}(X)\cong\Hilb^G(\A_k^3)$ for $G=N/N_0\subset\GL(3,k)$ (Remark \ref{Hilb FB iso}).
\end{rem}

\begin{rem}\label{only essential}
Conversely, \cite[Problem 1.5]{CMDY25} asks whether only essential divisors appear on $\FB_{(\infty)}(X)$.
For $X$ as in Theorem \ref{Hilb canonical}, by (ii) this asks whether a small resolution of $\Hilb^G(\A_k^3)$ is a Hilbert basis resolution in the sense of \cite{SS23}.
The answer is yes if $G\subset\SL(3,k)$, or if $G$ is of type $\frac1r(1,1,r-1)$, in which case $\Hilb^G(\A_k^3)$ is the Danilov resolution \cite{Ked11} (see \cite[Section 3]{SS23}).
The answer is no for the other groups of type (2), which can be seen by counting the cones of the fan \cite[Theorem 6.4]{Ked11}.
It is also no for type (5) \cite[Section 2]{SS23} and for the group $\frac14(1,2,3)$ of type (3) \cite[Example 6.16]{CMDY25}.
\end{rem}

\bibliographystyle{alpha}
\bibliography{Hilb}

\end{document}